\documentclass{rQUF2e}

\usepackage{epstopdf}
\usepackage{subfigure}
\usepackage{adjustbox}
\usepackage{graphicx}

\theoremstyle{plain}

\theoremstyle{definition}

\theoremstyle{remark}

\usepackage{comment}

\usepackage{todonotes}

\def\EE{\mathbb{E}}

\usepackage{tikz}
\usetikzlibrary{positioning,arrows,calc,fit}

\begin{document}


\title{{\textit{Learning a Functional Control \\for High-Frequency Finance}}}

\author{
     L. LEAL,\textsuperscript{\dag}
     M. LAURIERE,\textsuperscript{\dag}  
     C.-A. LEHALLE\textsuperscript{\ddag}
     \\
}
\affil{
     \textsuperscript{\dag}Department of Operations Research
 and Financial Engineering, Princeton University, Princeton, NJ - 08544 \\
 \{lleal,lauriere\}@princeton.edu\\
    \textsuperscript{\ddag} Capital Fund Management, Paris,
   and Imperial College, London. 
     23 rue de l'Université,
     75007 Paris, France\\
   charles-albert.lehalle@cfm.fr
}

\maketitle 


\begin{abstract}
We use a deep neural network to generate controllers for optimal trading on high frequency data. For the first time, a neural network learns the mapping between the preferences of the trader, i.e. risk aversion parameters, and the optimal controls. An important challenge in learning this mapping is that, in intra-day trading, trader's actions influence price dynamics in closed loop via the market impact. The exploration--exploitation tradeoff generated by the efficient execution is addressed by tuning the trader's preferences to ensure long enough trajectories are produced during the learning phase. The issue of scarcity of financial data is solved by transfer learning: the neural network is first trained on trajectories generated thanks to a Monte-Carlo scheme, leading to a good initialization before training on historical trajectories. Moreover, to answer to genuine requests of financial regulators on the explainability of machine learning generated controls, we project the obtained ``blackbox controls'' on the space usually spanned by the closed-form solution of the stylized optimal trading problem, leading to a transparent structure. For more realistic loss functions that have no closed-form solution, we show that the average distance between the generated controls and their explainable version remains small. This opens the door to the acceptance of ML-generated controls by financial regulators.
\end{abstract}


\begin{classcode}
C450; C580; C610; G170; C630 \end{classcode}

\clearpage

\section{Introduction}

Financial mathematics are using stochastic control to ensure that market participants are operating as intermediaries and not as unilateral risk takers: investment banks have to design risk replicating strategies, systemic banks have to ensure they have plans in case
of fire sales triggered by economic uncertainty, asset managers have to balance risk and returns in their portfolios and brokers have to be sure investors’ large buy or sale orders are executed without distorting the prices. The latter took a primary role in post-2008 markets since participants understood preservation of liquidity is of primary importance. 

Our paper addresses this last case, belonging to the academic field of optimal trading, initially introduced by \citep{ac01} and \citep{bl98}, and then extended in a lot of ways, from sophisticated stochastic control \citep{bdl11} to Gaussian-quadratic approximations allowing to obtain closed-form solutions like \citep{cj16} or \citep{cjp15}, or under a self-financing equation context, as in \citep{cw19}. Learning was introduced in this field either via on-line stochastic approximation \citep{llp13} or in the context of games with partial information \citep{cl18}. More recently, Reinforcement Learning (RL) approaches have been proposed, like in \citep{gm19} to face the high dimensionality of optimal trading of portfolios, or in \citep{ml19} to adapt the controls in an online way; see also~\citep{charpentier2020reinforcement} for an overview.

For optimal control problems in continuous time, the traditional approach starts with the cost function and derives, through dynamic programming, a Hamilton-Jacobi-Bellman (HJB) equation involving the optimal control and the value function (that is the optimum of the Q-function in RL).
This Partial Differential Equation (PDE) can be explicitly written only once the dynamics to be controlled are stylized enough. When it does not have a closed form solution, this PDE can be solved either by a deterministic method (like a finite difference scheme) that is limited by the dimension of the problem, or by RL approaches, like in \citep{bgtw19} for deep hedging or in \citep{gmp20, ttmz19} for optimal trading. In contrast with \citep{ttmz19} mentioned above, we focus less on the use of signals in algorithmic trading and more on the risk control dimension, going into neural network {\em explainability} aspects that, to the best of our knowledge, have not been explored before. 

In the broader literature, deep learning techniques for PDEs and Backward Stochastic Differential Equations (BSDEs) have recently attracted a lot of interest \citep{MR3874585,hure2019some,rpk19,whj17,hjw18} and found numerous applications such as price formation~\citep{sirignano2019universal}, option pricing~\citep{liu2019pricing} or financial systems with a continuum of agents~\citep{fouque2019deep,carmona2019convergence-I,carmona2019convergence-II}. 

In this paper, we propose to skip most of the previous steps: we directly learn the optimal control on the discrete version of the dynamics and of the cost function that are usually giving birth to the HJB, which we no longer have to derive. This approach, also used e.g. in~\citep{MR2137498,he16} for optimal control, gives us the freedom to address dynamics for which deriving explicitly a PDE is not possible, and to learn directly on trajectories coming from real data. In optimal trading, the control has a feedback loop with the dynamics via price impact: the faster you trade, the more you move the price in a detrimental way \citep{bill15}. In our setup the controller sets the trading speed of an algorithm in charge of buying or selling a large amount of stock shares, solving the trade-off between trading fast to get rid of the exposure to uncertainty of future prices and trading slow to not pay too much in market impact. 

For numerical applications we use high frequency (HF) data from the Toronto Exchange (TSX) on 19 stocks over two years, each of them generating around 6,500 trades (i.e. data points) per day. We compare different controls: controls generated by a widely known stylized version of the problem that has a closed form formula \citep{cjp15}, controls learned on simulated dynamics whose parameters are estimated on the data in different ways (with or without stationary assumptions), and controls learned on real data. For the latter case, we transfer the learning from simulated data to real data, to make sure we start with a good initialization of the controller.

In our setup, a controller maps the state space describing price dynamics and the trader's inventory to a trading speed that leads to an updated version of the state, and so on until the end of the trading day. The same controller is used at every decision step, which corresponds to one decision every 5 minutes. When the controller is a neural net, it is used multiple times before the loss function can be computed and improved via back propagation, in the spirit of \citep{la98}.

Our first main contribution is that we not only train a neural net for given values of the end-user’s preferences, but we also train a neural net having two more inputs that are the trader's preferences (i.e. risk aversion parameters), so that this neural net is learning the mapping between the preferences (i.e. hyper-parameters of the control setup) and the optimal controls.
To the best of our knowledge, it is the first time that a neural net is performing this kind of ``functional learning" of risk aversion parameters: the neural net learns the optimal control for a range of cost functions that are parametrized by the risk preferences. In the paper we call it a ``\emph{multi-preferences neural net}".

The second major contribution of this paper is the way we compare the generated controls in a meaningful way through {\em model distillation}. Model distillation and, more precisely, model translation methods seek to replicate the DNN behavior across an entire dataset (see Tree-based \cite{fh17,tchkg18,zymw19}, Graph-based \cite{zcwz16,zcswz17} and Rule-based \cite{ms17,hdr18} methods). Our paper introduces a model translation method focused on the idea of global approximation (see \cite{tchkg18} for a discussion on global additive explanations).

Our global distillation approach paves the way to methods which satisfy the request of financial regulators about the \textit{explainability of learned controls}. We start by the functional space of controls spanned by the closed-form solution of the stylized problem: they are non-linear in the remaining time to trade and affine in the remaining quantity to trade (see \citep{cl18} for a description of the relationship between the optimal controls and the space generated by the $h_1(t)$ and $h_2(t)$ defined later in the paper). Hence, we project the learned controls on an affine basis of functions, for each slice of remaining time to trade $T-t$. In doing so, we provide evidence of the distance between the black-box controls learned by a neural network, and the baseline control used by banks and regulators. We can compare the controls in this functional space and, thanks to the $R^2$ of the linear regressions, we know the average distance between this very transparent representation of the controls and the ones generated by the neural controllers. In practice, we show that when the loss function is mainly quadratic, the learned controls, even when they are trained on real data, are spanning roughly the same functional space as the closed-form ones. End-users and regulators can consider them with the same monitoring tools and analytics than more standard (but sub-optimal) controls. When the loss function is more realistic but not quadratic anymore, taking into account the mean reverting nature of intra-day price dynamics \citep{ll12}, the generated controls may or may not belong to the same functional space.

The structure of the paper is as follows: Section~\ref{sec2} presents the optimal execution problem setup, focusing on the loss function and on the closed-form formula associated with the state-of-the-art optimal execution problem. It also introduces the architecture of the neural networks we use and our learning strategies. Section~\ref{sec3} describes the dataset, stylized facts of intra-day price dynamics that should be taken into account, and motivates learning in a data-driven environment. Section~\ref{sec3} also presents the numerical results and our way to tackle explainability of the generated controls. Conclusions and perspectives are provided in Section~\ref{sec5}.

\section{The optimal execution model}
\label{sec2}

\subsection{Dynamics and payoff}
\label{sec2.1}

Optimal trading is dealing with an agent who would like to execute a large buy or sell order in the market before a time horizon $T$. Here, their control is a trading speed $\nu_t \in \mathbb{R}$. The center of the problem is to find the balance between trading too fast (and hence move the price a detrimental way and pay trading costs, as modelled here in \eqref{eq:1.1} and \eqref{eq:1.3}) and trading too slow (and being exposed to not finishing the order and be exposed to uncertainty of future prices, reflected in \eqref{eq:1.2} and \eqref{eq:1.4}). 
They maximize their wealth while taking into account a running cost and a final cost of holding inventory, and they are constrained by the dynamics of the system, which is described by the evolution of the price $S_t$, their inventory $Q_t$ (corresponding to a number of stocks), and their wealth $X_t$ (corresponding to a cash amount). 
In the market, the price process for the asset they want to trade evolves according to
\begin{equation}
\Delta S_t := S_{t+1} - S_{t} = \alpha_t \nu_t \Delta t +\sigma \sqrt{\Delta t} \, \epsilon_t, 
\label{eq:1.1}
\end{equation}
for $t=0,\dots,T-1$,  where $\alpha_t>0$ is a parameter that accounts for the drift, $\Delta t$ is the size of a time step, $\sigma>0$ is a constant volatility term, taken to be the historical volatility of the stock being traded. For transfer learning on historical data, the database contains real price increments and hence the realized volatility is (implicitly) of the same order of magnitude. $\epsilon \sim \mathcal{N}(0,1)$ is a noise term.

The state of the investor at time $t$ is described by the tuple $(T-t,Q_t,X_t)$. $Q_t$ is their inventory, or number of stocks, and $X_t$ is their wealth at time $t$, taken as a cash amount. To isolate the inventory execution problem from other portfolio considerations, the wealth of the agent at time $0$ is considered to be $0$, i.e. $X_0=0$. Suppose the agent is buying inventory (the selling problem is symmetric). Then, $Q_0<0$ and the trading speed $\nu$ will be mostly positive throughout the trajectory. The evolution of the inventory is given by:
\begin{equation}
\Delta Q_t := Q_{t+1} - Q_{t} =\nu_t \Delta t, 
\hspace{3mm} \text{for } t=0,\dots,T-1. \label{eq:1.2}
\end{equation}
The wealth of the investor evolves according to
\begin{equation}
\Delta X_t := X_{t+1} - X_{t} = -\nu_t(S_t+\kappa \cdot \nu_t)\Delta t, 
\label{eq:1.3}
\end{equation}
for $t=0,\dots,T-1$, where $\kappa>0$ is a constant, and the term $\kappa \cdot \nu_t$ represents the temporary market impact from trading at time $t$. As we can see, this is a linear function of $\nu_t$ that acts as an increment to the mid-price $S_t$. It can also be seen as the cost of ``crossing the spread'', or a transaction cost.

The cost function of the investor fits the  framework considered e.g. in~\citep{cjp15}. The agent seeks to maximize over the trading strategy $\nu$ the reward function given by:
\begin{equation}
J_{A,\phi}(\nu) :=  \EE \Big[ X_T + Q_T S_T-A |Q_T|^{\gamma} -\phi \sum_{t=0}^T |Q_t|^{\gamma} \Big], \label{eq:1.4}
\end{equation}
with $\gamma>0$. We mostly focus on the standard case $\gamma=2$, but we will also consider the case where $\gamma=3/2$, which better takes into account the sub-diffusive nature of intra-day price dynamics, see \citep{ll12}. $X_T, S_T$, and  $Q_T$ are the random variables parameterizing the terminal value of the stochastic processes described in the dynamics \eqref{eq:1.1}--\eqref{eq:1.3}. $A>0$ and $\phi>0$ are constants representing the risk aversion of the agent. $A$ penalizes holding inventory at the end of the time period, and $\phi$ penalizes holding inventory throughout the trading day. Together, they parametrize the control of the optimal execution model and stand for the agent's preferences.  

The trading agent wants to find the optimal control $\nu$ for the cost \eqref{eq:1.4} subject to the dynamics \eqref{eq:1.1}--\eqref{eq:1.3}. In the sequel, we solve this problem using a neural network approximation for the optimal control. As a benchmark for comparison, we will use state-of-the-art closed-form solution for the corresponding continuous time problem, which is derived through a PDE approach. The closed-form solution of the PDE is well-defined only when the data satisfies suitable assumptions.
The neural net learns the control $\nu_t$ directly from the data, while finding a closed-form solution for the PDE is not always possible.

It can be shown, see~\citep{cj16}, that the optimal control $\nu^*$ is obtained as a linear function of the inventory, which can be written explicitly as:
\begin{equation}
\nu^* (t,q) = \frac{h_1(t)}{2 \kappa}  + \frac{\alpha +h_2(t)}{2 \kappa}  q,
\label{eq:1.5}
\end{equation}
\noindent where $h_2$ and $h_1$ are the solutions of a system of Ordinary Differential Equations (ODEs).

\subsection{Non-parametric solution}
\label{sec2.2}

\paragraph{\textbf The neural network setup}
The neural network implementation we propose precludes all the usual derivations in the \citep{cj16} framework. We no longer need to find the PDE corresponding to the stochastic optimal control problem, we no longer need to break it down into ODEs and solve the ODE system, and we can try to directly approximate the optimal control. The deep neural network approximation looks for a control minimizing the agent's objective function \eqref{eq:1.4} while being constrained by the dynamics \eqref{eq:1.1}--\eqref{eq:1.3}, without any further derivations.

We define \textit{one single neural network $f_{\theta}(\cdot)$ to be trained for all the time steps}. Each iteration of the stochastic gradient descent (SGD) proceeds as follows. Starting from an initial point $(S_0,X_0,Q_0)$, we simulate a trajectory using the control $\nu_t = f_{\theta}(t, Q_t)$. Based on this trajectory, we compute the gradient of associated cost with respect to the neural network's parameters $\theta$. Finally, the network's parameters are updated based on this gradient. In our implementation, the learning rate is updated using Adaptive Moment Estimation (Adam)~\citep{kb14}, which is well suited for situations with a large amount of data, and also for non-stationary, noisy problems like the one under consideration. 

In the Monte-Carlo simulation mode, we generate random increments of the Brownian motion through the term $\sigma \sqrt{\Delta t} \, \epsilon$, where $\epsilon\sim \mathcal{N}(0,1)$ comes from the standard Gaussian distribution. Using the Gaussian increments, along with the $\nu_t$ obtained from the neural network, and the three process updates ($\Delta S_t,\Delta X_t,\Delta Q_t$), we can compute the new state variables in discrete time:
\begin{align}
S_{t+1} &= S_t + \alpha_t f_{\theta}(t,Q_t) \Delta t +\sigma \sqrt{\Delta t} \, \epsilon, \\
X_{t+1} &= X_{t}  - f_{\theta}(t,Q_t) (S_t+\kappa_t f_{\theta}(t,Q_t)) \Delta t,\\
Q_{t+1} &= Q_t + f_{\theta}(t,Q_t) \Delta t.
\end{align}
The new state variable $Q_{t+1}$, along with the time variable, is then going to serve as input to help us learn the control at the next time step. This cycle continues until we have reached the time horizon $T$. Since we are using the same neural network $f_\theta$ for all the time steps, we expect it to learn according to both time and inventory. Note that each SGD iteration requires $T$ steps because we need the neural net to have controlled the $T$ steps before being able to compute the value function for the full trajectory. It then backpropagates the error on the terminal value function across the weights of the same neural net we use at each of the $T$ iterations.

\paragraph{\textbf Neural network architecture and parameters} 
 The architecture consists of three dense hidden layers containing five nodes each, using the hyperbolic tangent as activation function. Each hidden layer has a dropout rate of $0.2$. We add one last layer, without activation function, that returns the output. The learning rate is $\eta=5e^{-4}$; the mini-batch size is 64; the tile\footnote{The tile size ($ts$) stands for how many samples of the inventory we combine with each sampled pair. For $ts=3$, each pair $(S_0,\{\Delta W_t\}_{t=0}^T)$ becomes three samples: $(Q_0^j,S_0,\{\Delta W_t\}_{t=0}^T)$, for $j\in\{1,2,3\}$, where $j$ is the tile index. This is useful when using real data, because we can generate more scenarios than we would ordinarily be able to.} size is 3; and the number of SDG iterations is 100,000. Every 100 SGD iterations, we perform a validation step in order to check the generalization error, instead of evaluating just the mini-batch error. We used Tensorflow in our implementation.

Although the neural network's basic architecture is not very deep, each SGD step involves $T$ loops over the same neural network. This is because of the closed loop which the problem entails.  The number of layers is thus artificially multiplied by the number of time steps considered.

For the inputs, we have two setups. In the basic setup (as discussed above), the input is the pair $(t,Q_t)$. In the second case, which we call "\emph{multi-preferences neural network}", the input is the tuple $(t,Q_t,A,\phi)$ and the neural netowrk learns to minimize the $J_{A,\phi}(\cdot)$ cost functions for all $(A,\phi)$.
Each time we need to solve a given system, we set $A$ and $\phi$ to the desired value in the multi-preferences network and we do not need to relearn anything to obtain the optimal controls. For both neural nets, the output is the speed of trade $\nu_t$ for each time step $t$, intercalating with the variable updates, thus learning a controller for the full length of the trading day.

In order to learn from historical data, we first train the network on data simulated by  Monte Carlo and then perform transfer learning on real data.

\subsection{Explainability of the learned controls using projection on a well-known manifold}
\label{sec2.3}

We would like to compare the control obtained using the PDE solution with the control obtained using the neural network approximation. As stated in equation \eqref{eq:1.5}, the optimal control can be written as an affine function of the inventory $q_t$ at any point in time, for $t \in [0,T]$. 
The shape of this closed-form optimal control belongs to the manifold of functions of $t$ and $q$ spanned by $[0,T] \times \mathbb{R} \ni (t, q) \mapsto  h_1(t)/(2 \kappa)+(\alpha + h_2(t))/(2 \kappa)\cdot q \in \mathbb{R}$, where the $h_i(t)$ are non-linear smooth functions of $t$. To provide \emph{explainability of our neural controller} we project its effective controls on this manifold.
We obtain two non-linear functions ${\tilde h}_1(t)$ an ${\tilde h}_2(t)$ and a $R^2(t)$ measuring the distance between the effective control and the projected one at each time step $t$.

\paragraph{Procedure of projection on the "closed-form manifold"}
For each $t$, we form a database of all the learned controls $\nu_t$ mapped by the neural net to the remaining quantity $q_t$. 
It enables us to project this $\nu_t$ on $q_t$ using an Ordinary Least Squares (OLS) regression whose loss function is given by: 
\begin{equation}
    L(\beta_1(t),\beta_2(t)) = ||\nu(q(t)) - (\beta_1(t) + \beta_2(t)q(t))||^2_2,
\end{equation}
\noindent and do a global translation of our neural network control. 

The coefficients $\beta_{1}(t)$ and $\beta_{2}(t)$ of this regression $\nu_t = \beta_{1}(t) + \beta_{2}(t) q_t + \varepsilon_t$ can be easily inverted to give
\begin{equation}
    {\tilde h}_1(t) := 2\kappa\, \beta_1(t),\quad%
    {\tilde h}_2(t) := 2\kappa\, \beta_2(t) - \alpha
\end{equation}
for each $t \in [0,T]$. The $R^2(t)$ of each projection associated to $t$ quantifies the distance between the effective neural "black box" control and the explained one ${\tilde\nu}_t:={\tilde h}_1(t)/(2 \kappa)+(\alpha +{\tilde h}_2(t))/(2 \kappa)\cdot q_t$.
The curve of $R^2(t)$ represents how much of the non-linear functional can be projected onto the manifold of closed-form controls at each $t$. In practice we perform $T = 77$ OLS projections to obtain the curves of $h_1, h_2$ and $R^2$  (see Figure \ref{fig:2}) which provide explainability of our learned controls.

\paragraph{Users' preferences and exploration-exploitation issues}

The rate at which the agent executes the whole inventory strongly depends on the agent's preferences $A$ and $\phi$. 
When they are both large, the optimal control tends to finish to trade earlier than $T$, the end of the trading day.
Because of that, the trajectories of $Q_t$ approach zero in just a few time steps and after that large values of $Q_t$ are not visited during the training phase. As a consequence, in such a regime, it will be very difficult for the neural net to learn anything around the end of the time interval since it will no longer \emph{explore} its interactions with price dynamics to be able to emulate a better control.
For the case of the multi-preferences controller with 4 inputs, including $A$ and $\phi$, taking $(A,\phi)$ in a certain range ensures that the neural net will witness long enough trajectories to \emph{exploit} them. To this end, we took profit of the closed-form solution of the stylized dynamics and scanned the duration of trajectories for each pair $(A,\phi)$. 
In order to avoid the exploration - exploitation trade-off, we restricted the domain to values of $\phi$ smaller or equal to $0.007$, and values of A smaller or equal to $0.01$. This enabled the regression to learn the functional produced by the neural network. There will be some inventory left to execute at the end of the trading day for these parameters. For the pair $(A,\phi) = (0.01,0.007)$, we have less than $1\%$ of the initial inventory left to execute, yet it is enough for us to estimate the regression accurately.

\section{Optimal execution with data}
\label{sec3}

\subsection{Data description}
\label{sec3.1}

We use ticker data for trades and quotes for equities traded on the Toronto Stock Exchange for the period that ranges from Jan/2008 until Dec/2009, which results in 503 trading days. Both trades and quotes are available for Level I (no Level II available) for 19 stocks. They have a diverse range of industries, and daily traded volume. Our method can be directly used for all the stocks, either individually or simultaneously.\footnote{
For the market impact parameters $(\alpha, \kappa)$, we estimate the same value for all the stocks through a normalization, then we de-normalize the parameters using the seasonality profile of the stock used in the model.} Table \ref{tab:2.1} summarizes the database information:

\begin{table}
\begin{center}
\begin{minipage}{80mm}
\tbl{Descriptive statistics of the TSX stocks used.}
{\begin{tabular}{@{}rl}\toprule
  \multicolumn{1}{r}{Average \# Trades per Day:} & \multicolumn{1}{l}{124,988.4}      \\
\colrule
  \multicolumn{1}{r}{Average \# Trades per Stock:} & \multicolumn{1}{l}{3,308,902.2}     \\ 
  \colrule
\multicolumn{1}{r}{Permanent Market Impact $\alpha$:} & \multicolumn{1}{l}{$0.16 \cdot \frac{\text{avg. bin spread}}{\text{avg. bin volume}}\cdot \frac{1}{dt}$}     \\ \colrule

\multicolumn{1}{r}{Temporary Market Impact $\kappa$:} & \multicolumn{1}{l}{$0.24 \cdot \frac{\text{avg. bin spread}}{\text{avg. bin volume}}\cdot \frac{1}{dt}$}      \\ \colrule

\multicolumn{1}{r}{\% of 5-min Bins w/ Data:} & \multicolumn{1}{l}{91.5\%}    \\
\botrule
\end{tabular}}
\label{tab:2.1}
\end{minipage}
\end{center}
\end{table}
 
We partition both trade and quote files into smaller pickled files organized by date. We merge them day by day to calculate other interesting variables, such as: trade sign, volume weighted average price (VWAP), and the bid-ask spread at time of trade. We drop odd-lots from the dataset, as they in fact belong to a separate order book. Moreover, we only keep the quotes that precede a trade, or a sequence of trades. It saves us memory usage and processing times. 

The trade data, including the trade sign variable, is further processed into five minute bins. We have picked five minute intervals because, while the interval could be changed, here we want to make sure the data is not too sparse, and to provide additional risk control to the agent over a time horizon of one full day. TAQ data is extremely asynchronous and a situation might happen when there is not enough data in a given time interval. To avoid this problem, we aggregate the data into bins and focus our analysis on liquid stocks. This bin size is big enough for us to have data at each bin, and small enough to allow the agent to adjust the trading speed several times throughout the day. Since the market is open from 9:30 until 16:00, five-minute bin intervals, we have 78 time steps, which gives us 77 bins for which to estimate the control on each trading day. We have 91.5\% of bins containing data (see Table \ref{tab:2.1}). 

\subsection{Important aspects of financial data}

The stylized dynamics allowing a PDE formulation of the problem and a closed-form solution of the control make a lot of assumptions on the randomness of price increments, market impact and trading costs. Typically they assume independent and normally distributed, stationary, non-correlated returns, with no seasonality. However, these properties typically do not hold true for financial time series. See \citep{bbdg18,cgw93} for some interesting discussions.

\begin{figure}
\centering
\begin{minipage}{.45\linewidth}
  \includegraphics[width=\linewidth]{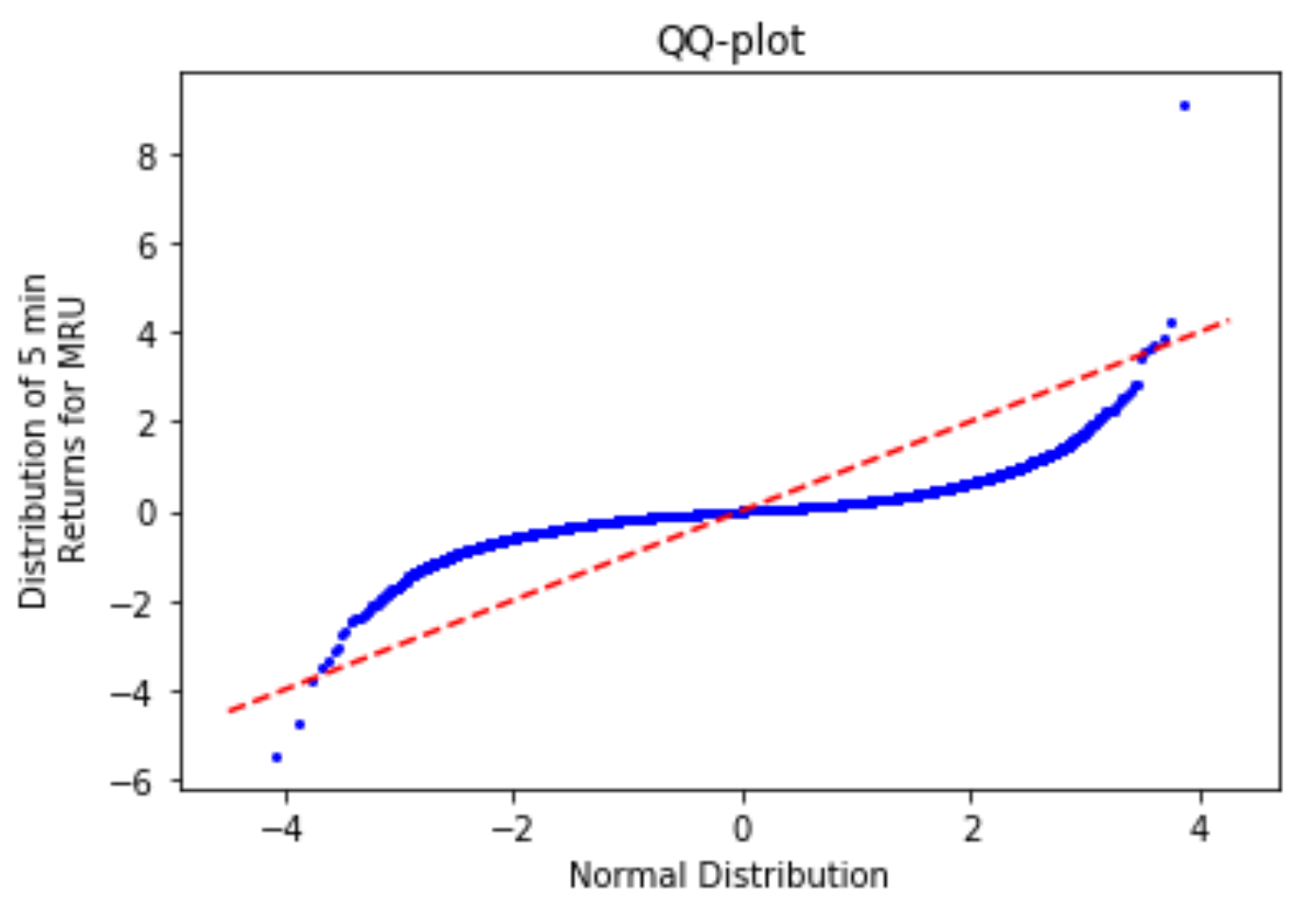}
\end{minipage}
\begin{minipage}{.44\linewidth}
  \includegraphics[width=\linewidth]{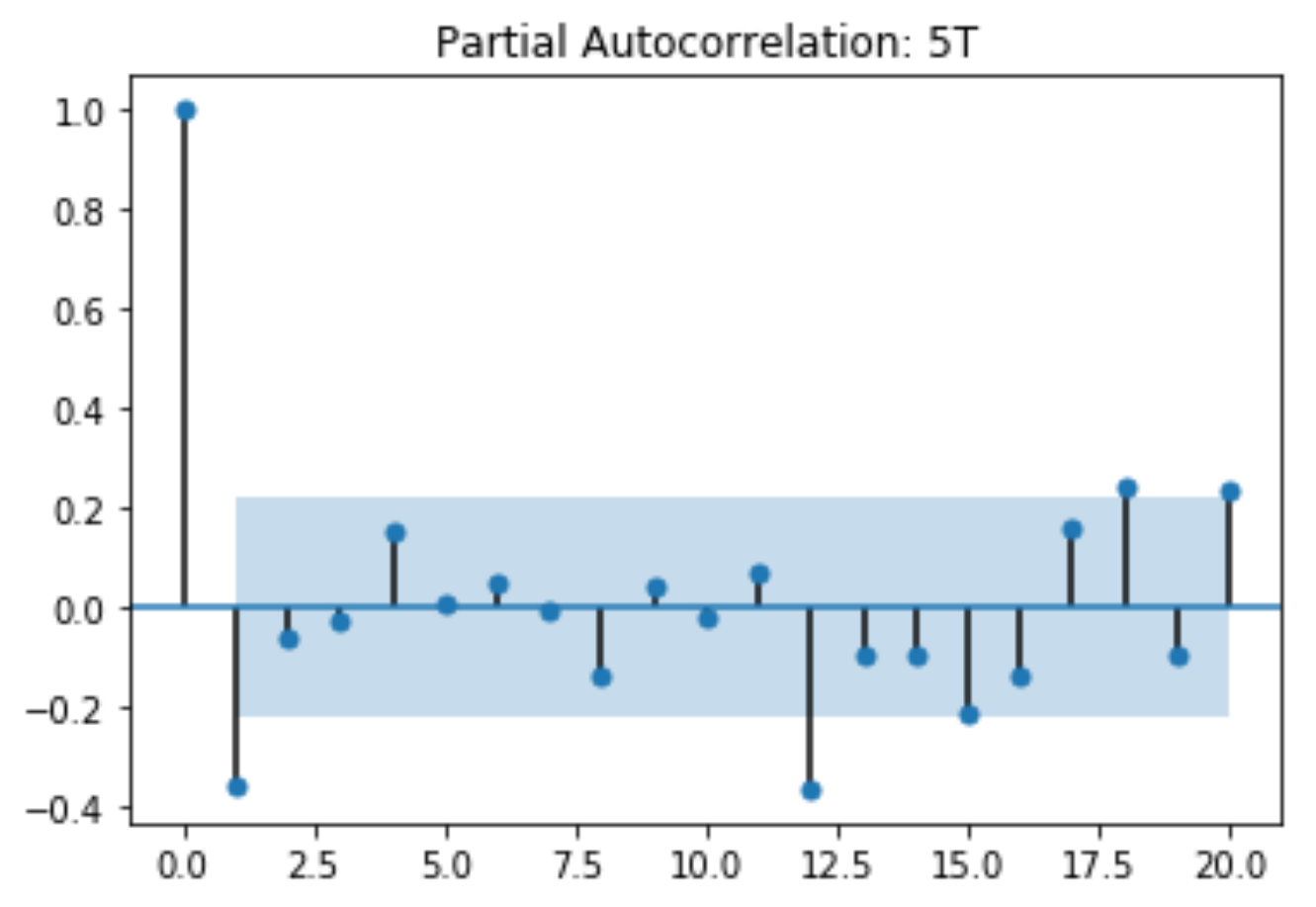}
\end{minipage}
\begin{minipage}{.53\linewidth}
  \includegraphics[width=\linewidth]{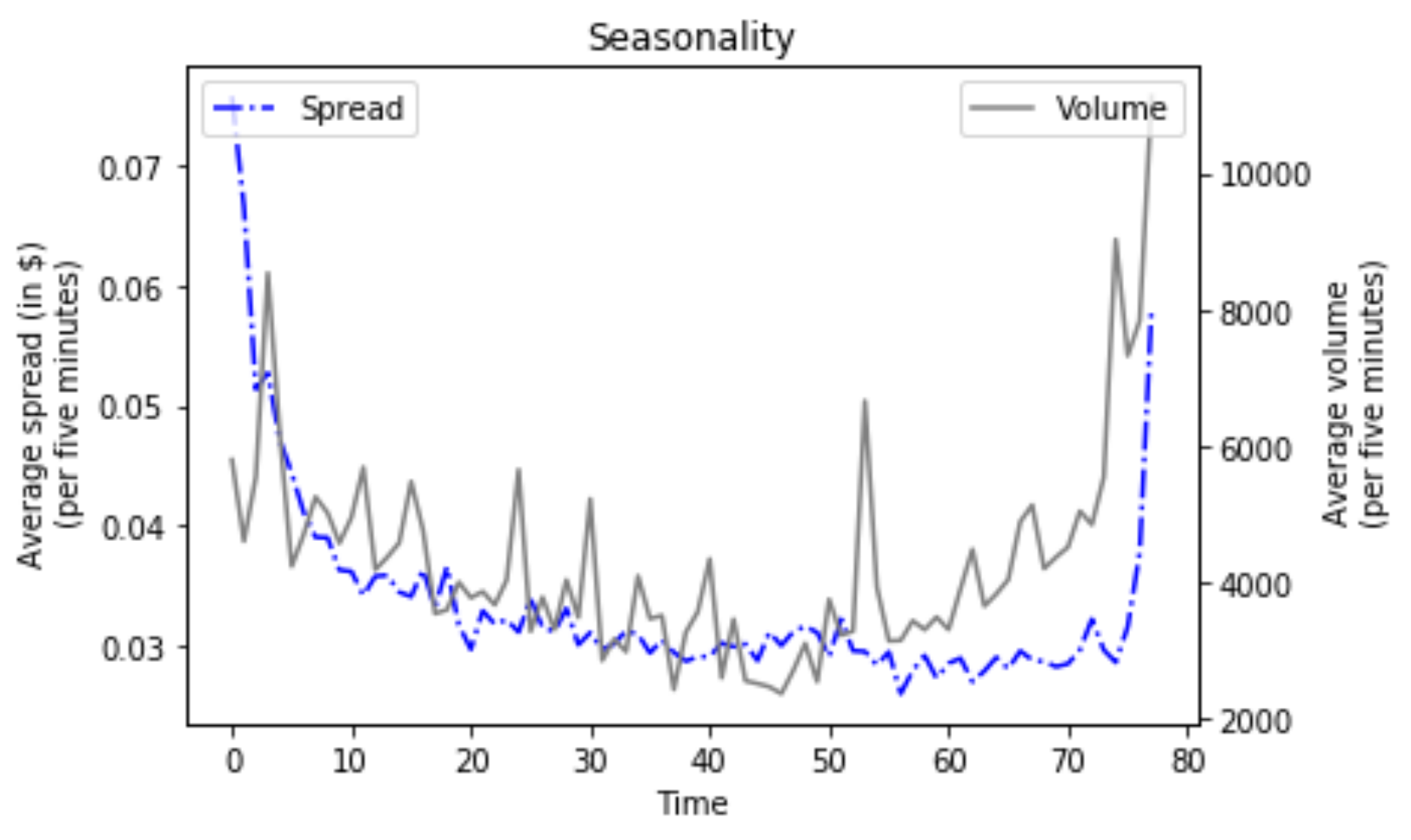}
\end{minipage}
\caption{QQ-plot (top left), auto-correlation function (top-right), intra-day seasonalities (bottom) of price dynamics}
\label{fig:1}
\end{figure}

Figure \ref{fig:1} shows the specific characteristics of financial data for the stock MRU. The left plot is a qq-plot, which tells us that the stock returns have a heavy tailed distribution. The same is confirmed for all stocks: we observe high kurtosis (ranging from $163$ to $11302$) for all the distributions of stock returns (which is a common metric for saying that a distribution has heavy tails with respect to a Gaussian distribution). The middle plot presents the five-minute lag auto-correlation profile. For this stock, we can see a mean-reversion tendency for five minutes and one hour. The first lag of the five-minute auto-correlation (ranging from $-0.017$ to $-0.118$) indicates all stocks present this feature. The right plot shows intra-day seasonality profiles: the bid-ask spread (blue, dashed curve), and the intra-day volume curve (full, gray curve).

As mentioned in section \ref{sec2.3}, the first step is to compare the output of our neural network model to the PDE solution, both using Monte Carlo simulated data. The second step is to understand how seasonality in the data might affect the training and the output. Finally, we continue our learning process on real data, in order to learn its peculiarities in terms of heavy tails, seasonality and auto-correlation.

\subsection{PDE \textit{vs} DNNs: from Monte Carlo to real data}
\label{sec4.1.1}

\begin{figure}
\centering
\begin{minipage}{.49\linewidth}
  \includegraphics[width=\linewidth]{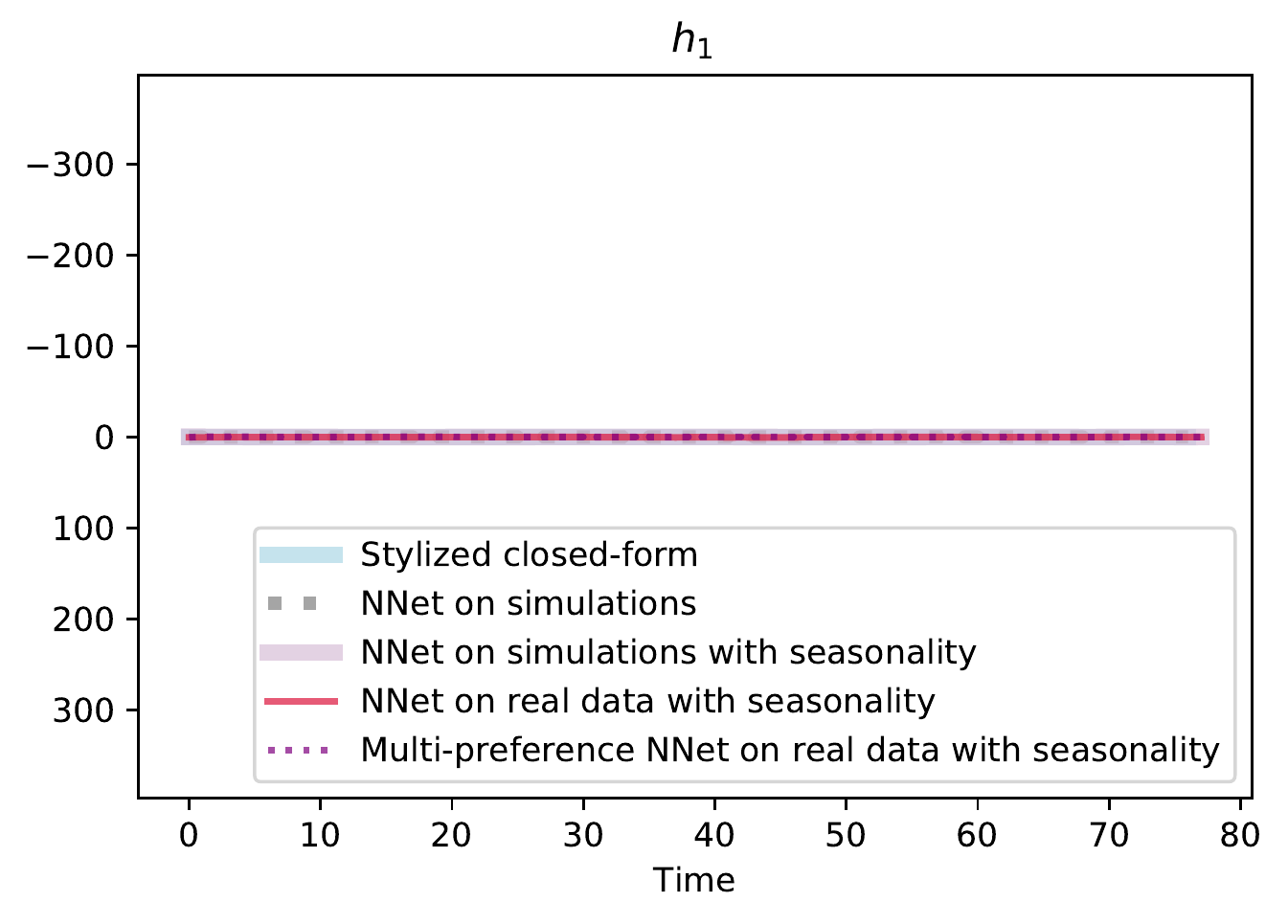}
\end{minipage}
\begin{minipage}{.49\linewidth}
  \includegraphics[width=\linewidth]{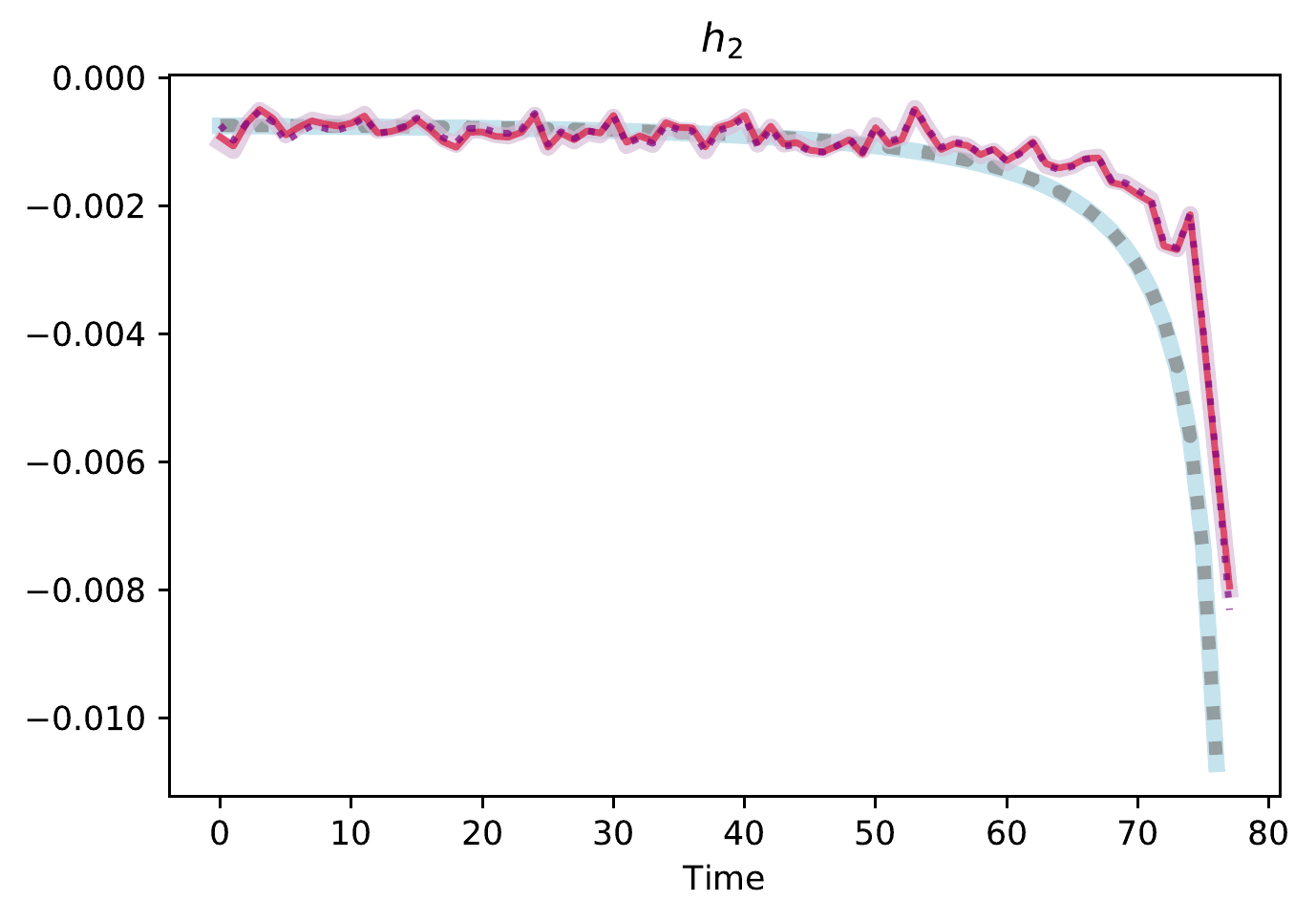}
\end{minipage}
\begin{minipage}{.49\linewidth}
  \includegraphics[width=\linewidth]{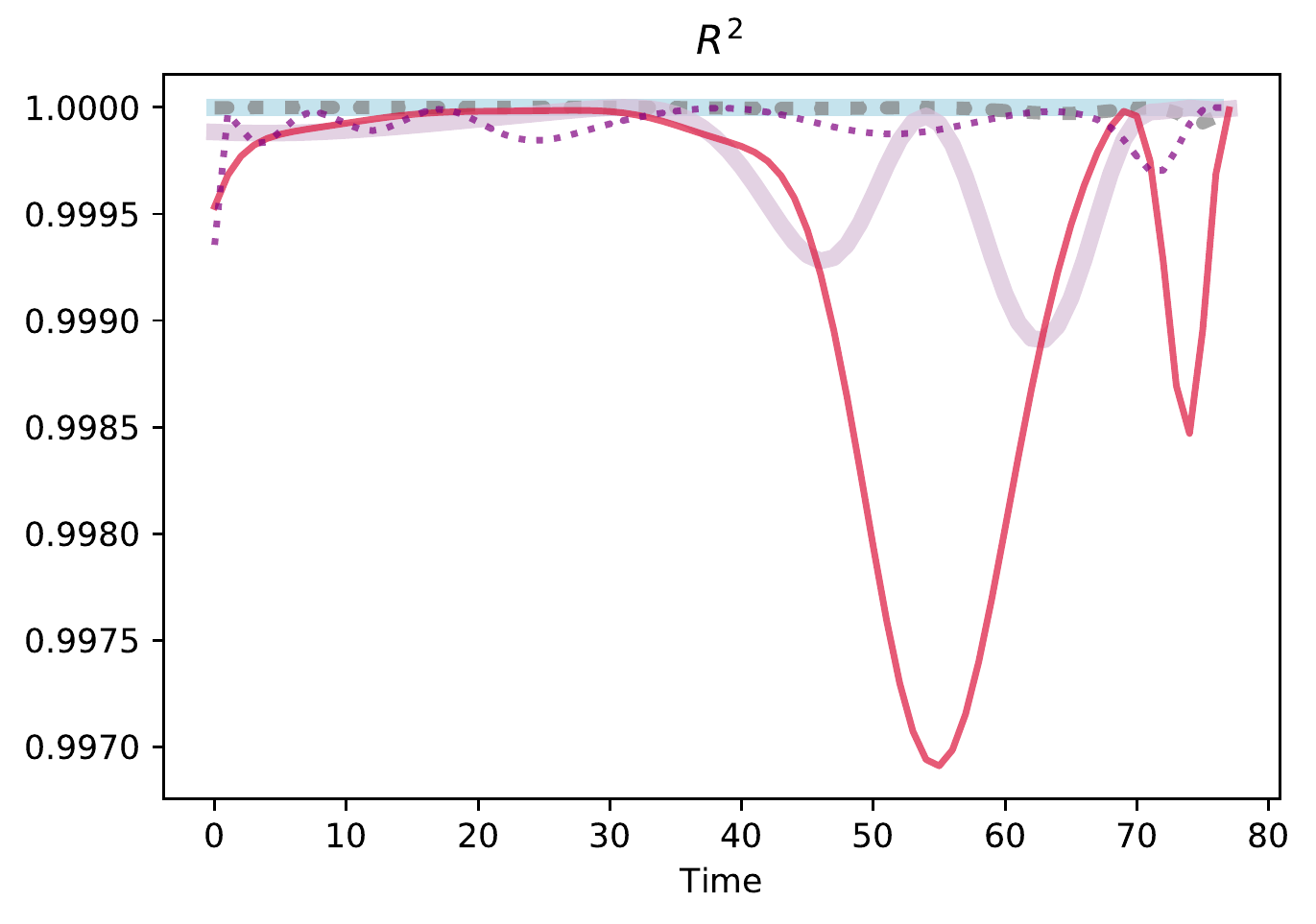}
\end{minipage}
\begin{minipage}{.49\linewidth}
  \includegraphics[width=\linewidth]{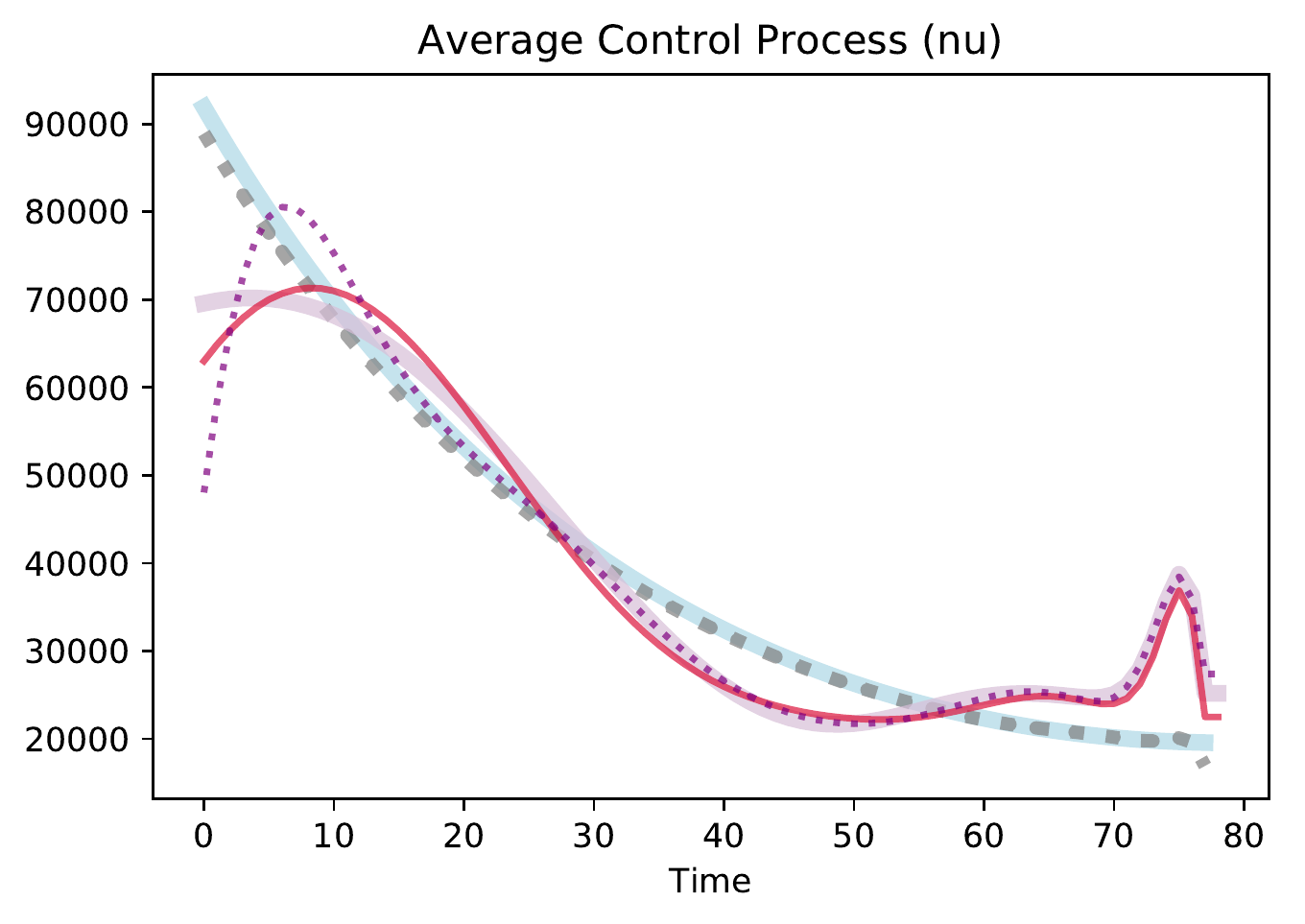}
\end{minipage}
\caption{Explainable parts of the control: $h_1(t)$, $h_2(t)$; fraction of the explained control $R^2(t)$, and average control given it is positive $\mathbb{E}(\nu(t)|\nu(t)>0)$ for different controllers.}
\label{fig:2}
\end{figure}

The results and our benchmark are summarized in Figure \ref{fig:2}. As stated in equation \eqref{eq:1.5}, the optimal control for the stylized dynamics and the PDE is linear in the inventory, hence the associated $R^2$ is obviously 1 for any $t$.
For the different neural nets, we perform the projections on the $(h_1,h_2)$ manifold and keep track of the $R^2$ curve.

During the learning, it is important that the neural network can observe full trajectories. When the risk aversion parameters are very restrictive, the closed-form solution and neural net trades so fast that the order is fully executed (i.e. the control stops) far before $t=T$. Because of that it is impossible to learn after this stopping time corresponding to $q_t=0$. Our workaround for this exploration - exploitation issue has been to use the closed-form solution to select a range of $(A,\phi)$ allowing the neural net to generate enough trajectories to observe the dynamics up to $t=T$. In Figure \ref{fig:2}, we use the pair $(A,\phi)=(0.01,0.007)$ with $\gamma=2$. We stress that, once the model has been trained on a variety of pairs $(A,\phi)$, we can then use it with new pairs of parameters.

In this Figure \ref{fig:2}, we mix closed-form controls on the stylized dynamics, neural controls on the same dynamics (using Monte Carlo simulations), neural controls in more realistic simulations (with an intra-day seasonality) and ultimately on real data. In particular, we use stock data (via transfer learning, see the thick purple and solid pink lines), and in order to make the model as close as possible to reality in a version that is in-between synthetic data and transfer learning (dotted pink curve). We include intra-day seasonality for both volume and bid-ask spread, and we emphasize how it is reflected in the learned control. We also superpose results for neural networks trained only on this regime of preferences and the controls of the multi-preference neural network, that is trained only once and then generates controls for any pair $(A,\phi)$ of preferences. The top panels represent the two components of the projection of the control on the "closed-form manifold", to enable comparison. The plot on the upper left of Figure \ref{fig:2} corresponds to the $h_1$ component. The $y$-axis scale is in the order of the median of the inventory path, to allow for a fair comparison to $h_2 \cdot Q$. It is naturally very close to zero, since there is no mean-field component in the market impact term. The upper right plot corresponds to the $h_2$ component. We notice that the stylized closed-form and the neural network simulation using Monte Carlo generated paths match perfectly, thus setting a great benchmark for the neural network. The improvements we add to the neural network model, simulations including seasonality, real data with seasonality and multi-preferences all overlap, indicating that the intra-day seasonality component is the most relevant when explaining differences in optimal execution speed when departing from the closed-form stylized model. The $R^2(t)$ curves (bottom-left panel) are flat at 1 for the closed form formula (since 100\% of the control belongs to this manifold), whereas it can vary for the other controls (see Figure \ref{fig:6}, and discussion in section \ref{sec:4.2}). Thanks to this projection, it is straightforward to compare the behaviour of all the controls. The $R^2(t)$ curves provide evidence that in this regime the neural controls are almost 100\% explained by the closed-form manifold (notice that the $y$-axis has a very tight range close to 1). It does not say that they are similar to the closed-form solutions of the stylized problem, but that they are linear in $q_t$ (but not in $t$) for each time step $t\in[0,T]$.

Saving the best for last, we present the optimal control for each different configuration in the bottom left plot of Figure \ref{fig:2}. We again stress the similarity between the neural network model on simulations versus the closed-form benchmark. We further observe the differences between the optimal controls for the improved setups. When we start taking into account the intra-day seasonality, the execution speed adjusts according to the volume traded and to the bid ask spread. It results in less trading when the spreads are large, and in more trading when there is more volumes. In adjusting to market conditions, we are able to improve on the existing execution benchmark.

\paragraph{Main differences between the learned controls}
Going from stylized dynamics to simulations with intra-day seasonality, and then to real price dynamics does not significantly change the multiplicative term of $q_t$, but it does shift the trading speed. This has been already observed in the context of game theoretical frameworks \citep{cl18}. When the seasonality is added in the simulation (thick purple curve) the deformation of the control exhibits oscillations that most probably stem from the seasonalities of the bottom panel of Figure \ref{fig:1}. The shift learned on real price dynamics (solid pink) amplifies some oscillations that are now less smooth; this is probably due to autocorrelations of price innovations. Moreover Figure \ref{fig:2} shows that the "functional learning" worked: the mono-preference neural network trained only on this $(A,\phi)$ and the multi-preferences one have similar curves.

\subsection{Nonlinear aspect of DNN controllers}
\label{sec:4.2}

\begin{figure}
\centering
\begin{minipage}{.48\linewidth}
  \includegraphics[width=\linewidth]{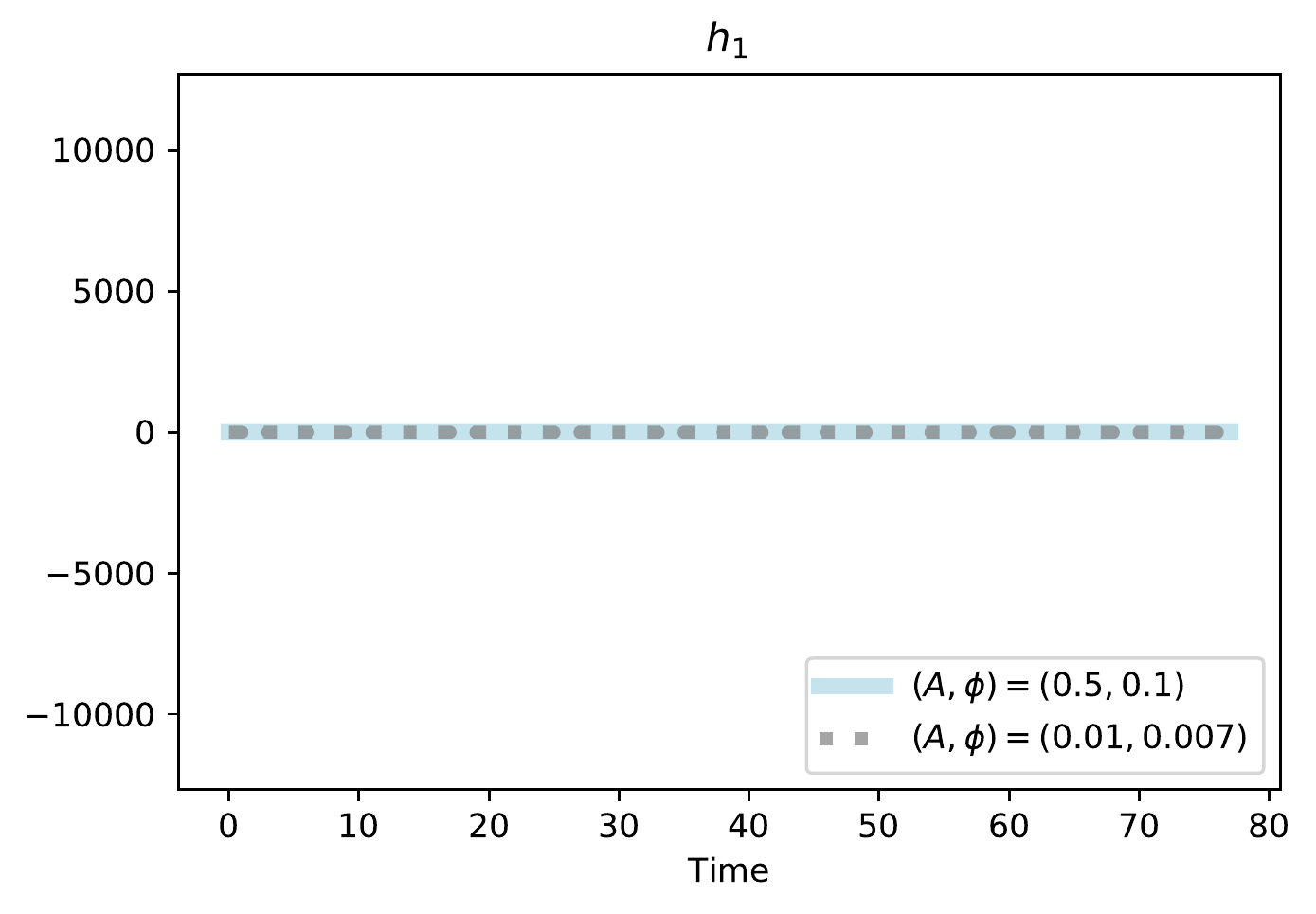}
\end{minipage}
\begin{minipage}{.48\linewidth}
  \includegraphics[width=\linewidth]{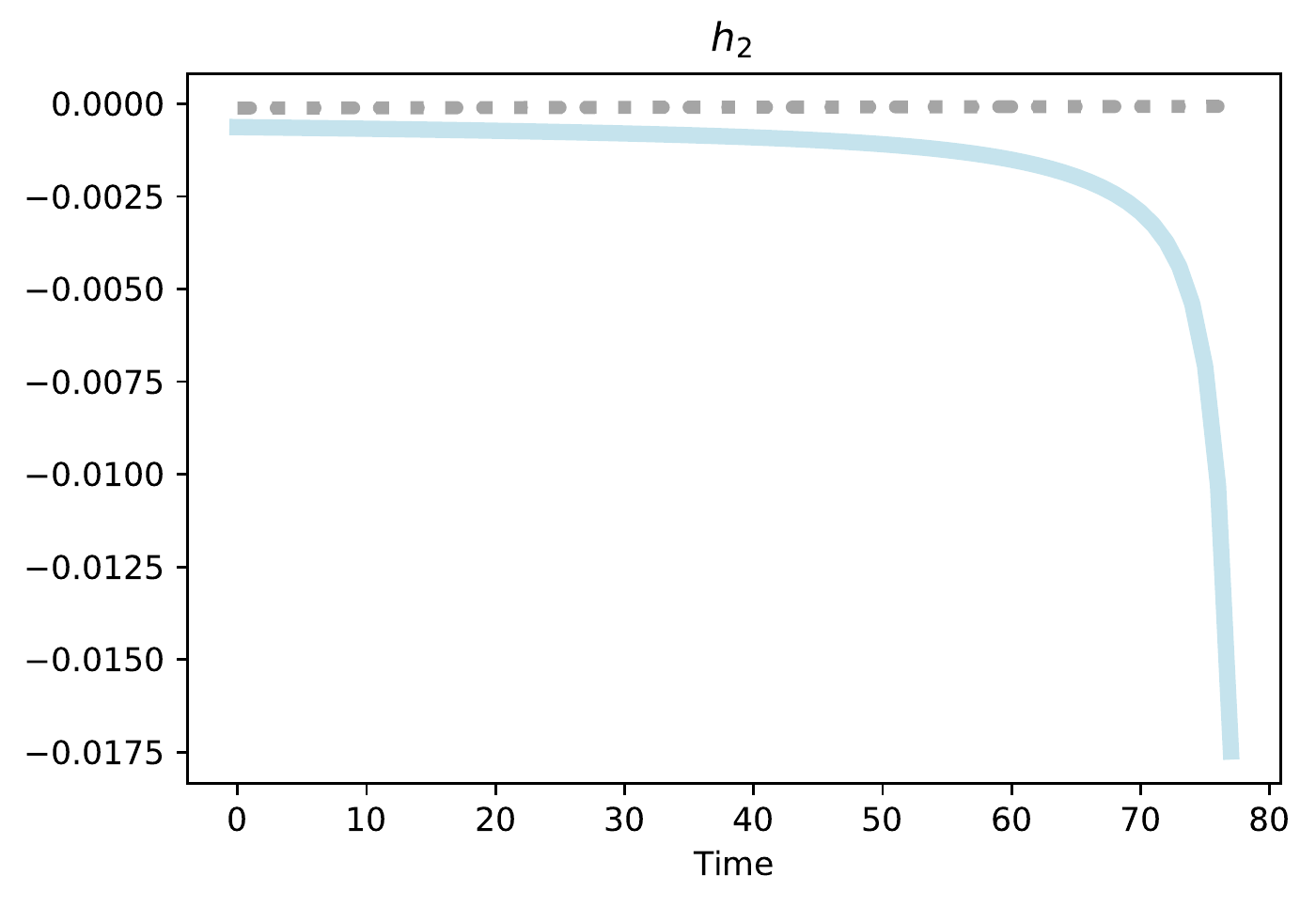}
\end{minipage}
\begin{minipage}{.48\linewidth}
  \includegraphics[width=\linewidth]{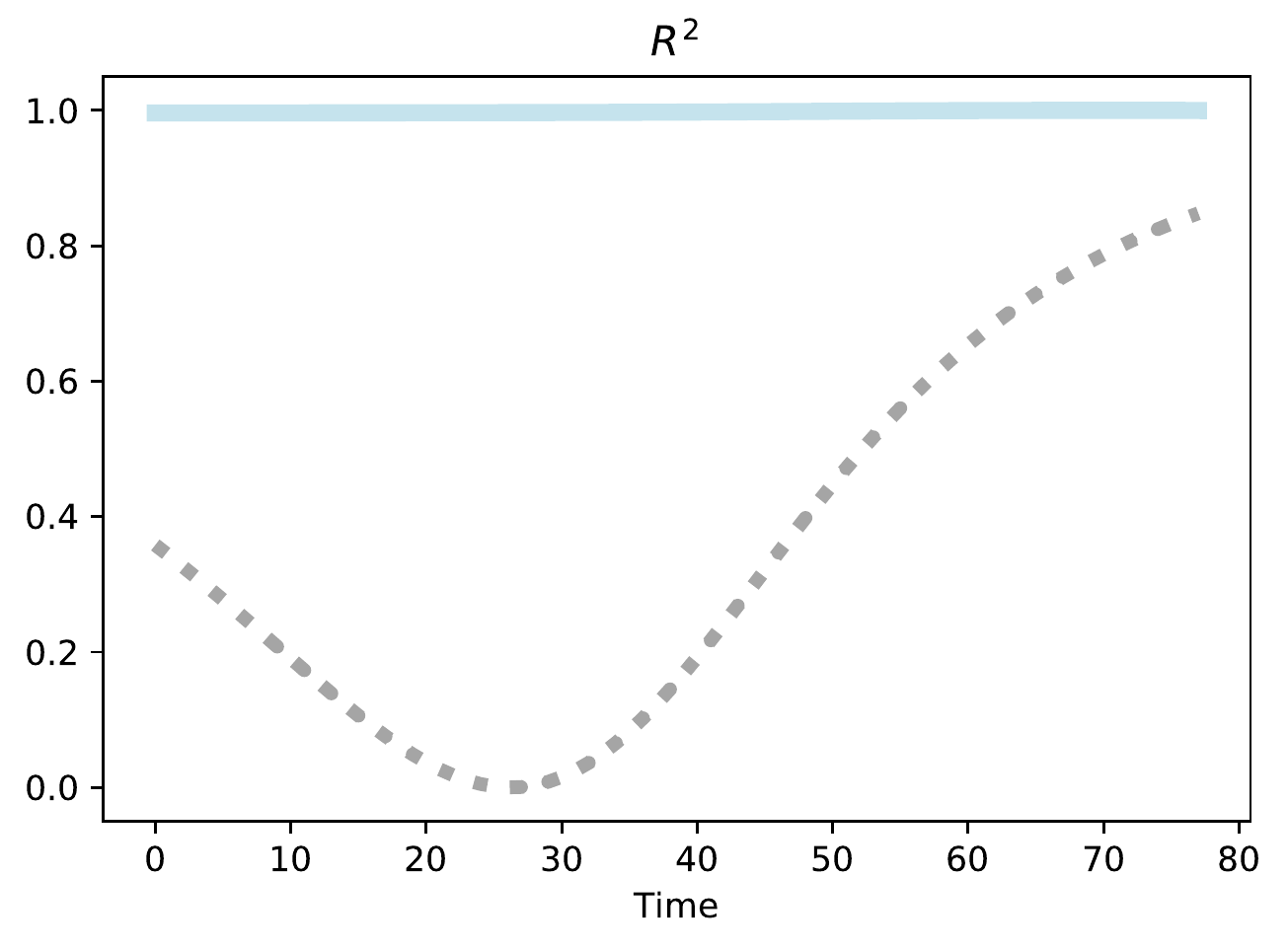}
\end{minipage}
\begin{minipage}{.48\linewidth}
  \includegraphics[width=\linewidth]{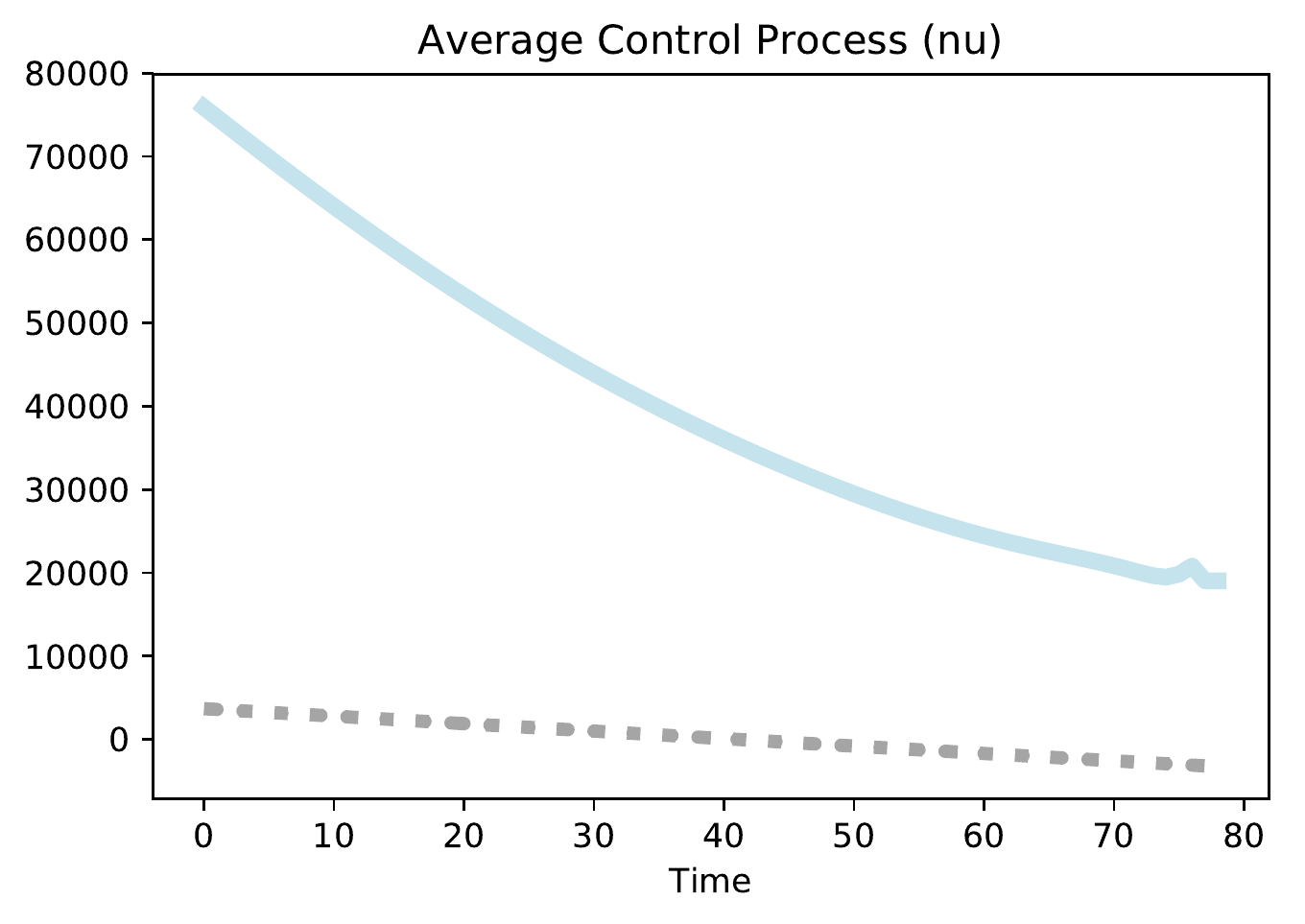}
\end{minipage}
\caption{$h_1$, $h_2$, $R^2$, $\nu$: for the case $\gamma=3/2$}
\label{fig:6}
\end{figure}

For non quadratic costs (we took $\gamma=3/2$ since it is a realistic case, see \cite{ll12}) we take $(A,\phi)=(0.5,0.1)$ to be in a comparable regime as the controls in Figure \ref{fig:2}, and $(A,\phi)=(0.01,0.007)$ which results in a non-linear control. The solid blue $R^2(t)$ curve in the bottom left of Figure \ref{fig:6} being again close to one shows that for this combination of parameters, the learned control remains extremely similar to the closed-form manifold. This kind of evidence can allow risk departments of banks and regulators to get more comfortable with this learned control and it can give guidelines to large institutional investors, such as pension funds, to execute their orders in the market with the view of controlling risk in an explicit fashion.

When $\gamma=3/2$, the chosen set of risk parameters has more influence in the learned control, which can now either be linear or non-linear. For the pair $(A,\phi)=(0.01,0.007)$, the absolute value of “proportional term” $h_2$ is very close to zero when using the same scale as the $h_2$ when $(A,\phi)=(0.5,0.1)$. Nevertheless, it is compensated by the part of the control that is outside of this space, as indicated by the $R^2$: at the start of the trading process 40\% of the control is explained by the “closed-form manifold”; it then decreases to 0\%, meaning that the learned control takes part fully out of this manifold; but it then increases back to 80\% at the end of the process. It seems that the best way to finish the trading is already well captured by the closed-form manifold. As we observe from the control plot in the bottom-right of Figure \ref{fig:6}, the parameter pair $(A,\phi)=(0.01,0.007)$ does not enforce trading speeds as fast as the pair $(A,\phi)=(0.5,0.1)$ when the loss function is sub-diffusive. Keeping the same preferences as before is thus making the trader less aggressive in a sub-diffusive environment. In order to have the same behavior in terms of control, they would need to behave in a more risk averse manner.

\subsection{Comparison to State-of-the-Art}

\begin{figure}
\centering
\begin{minipage}{.49\linewidth}
  \includegraphics[width=\linewidth]{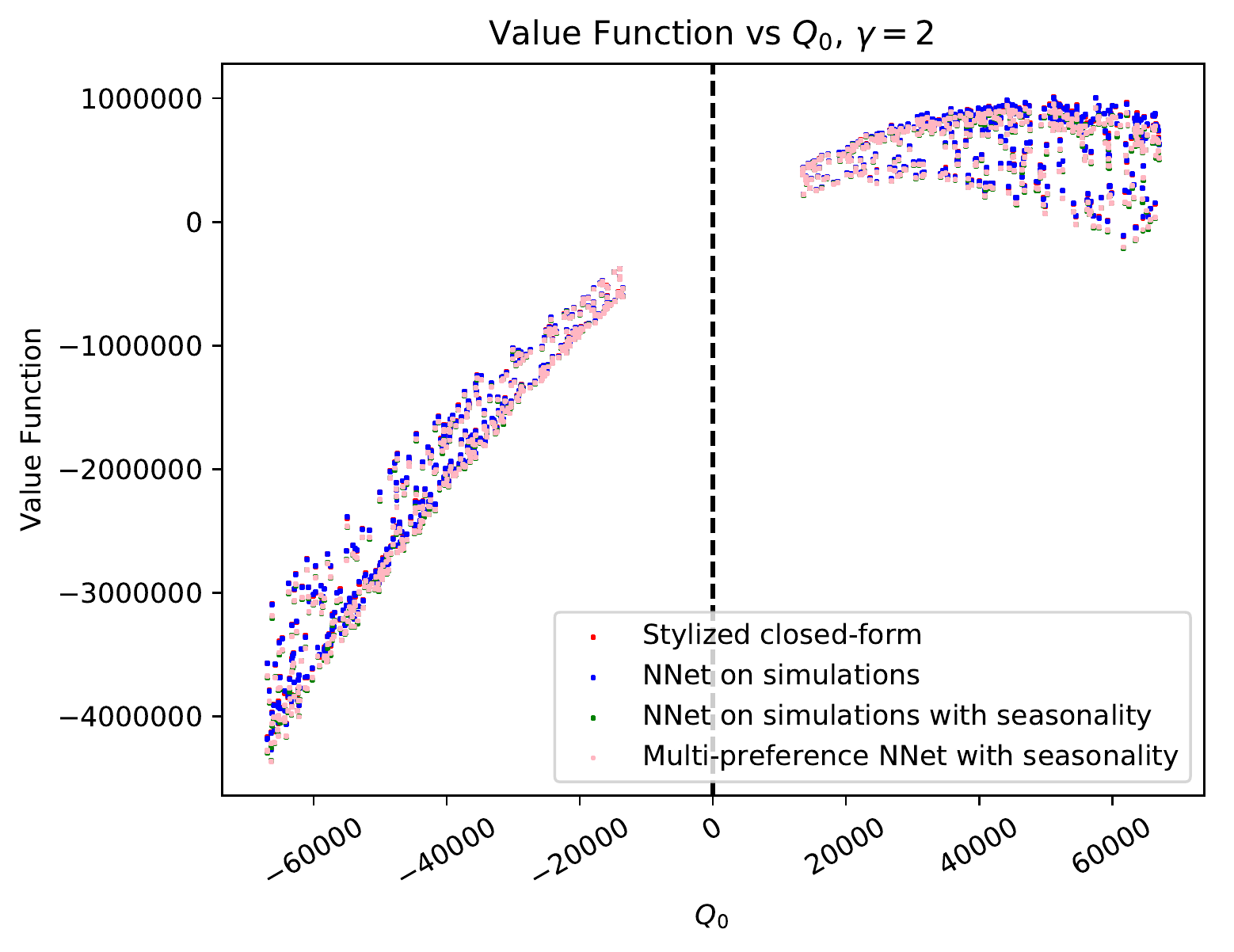}
\end{minipage}
\begin{minipage}{.49\linewidth}
  \includegraphics[width=\linewidth]{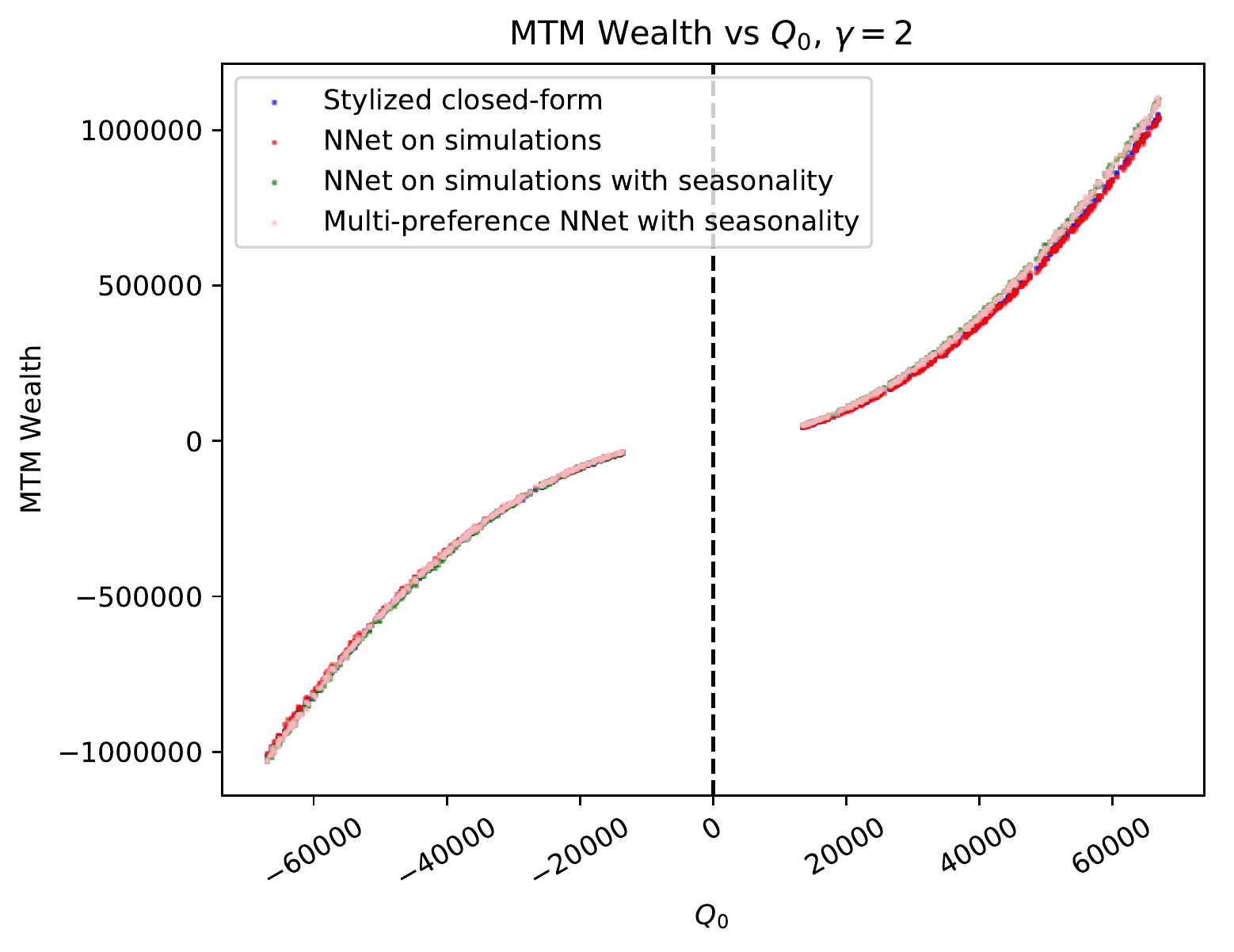}
\end{minipage}
\begin{minipage}{.49\linewidth}
  \includegraphics[width=\linewidth]{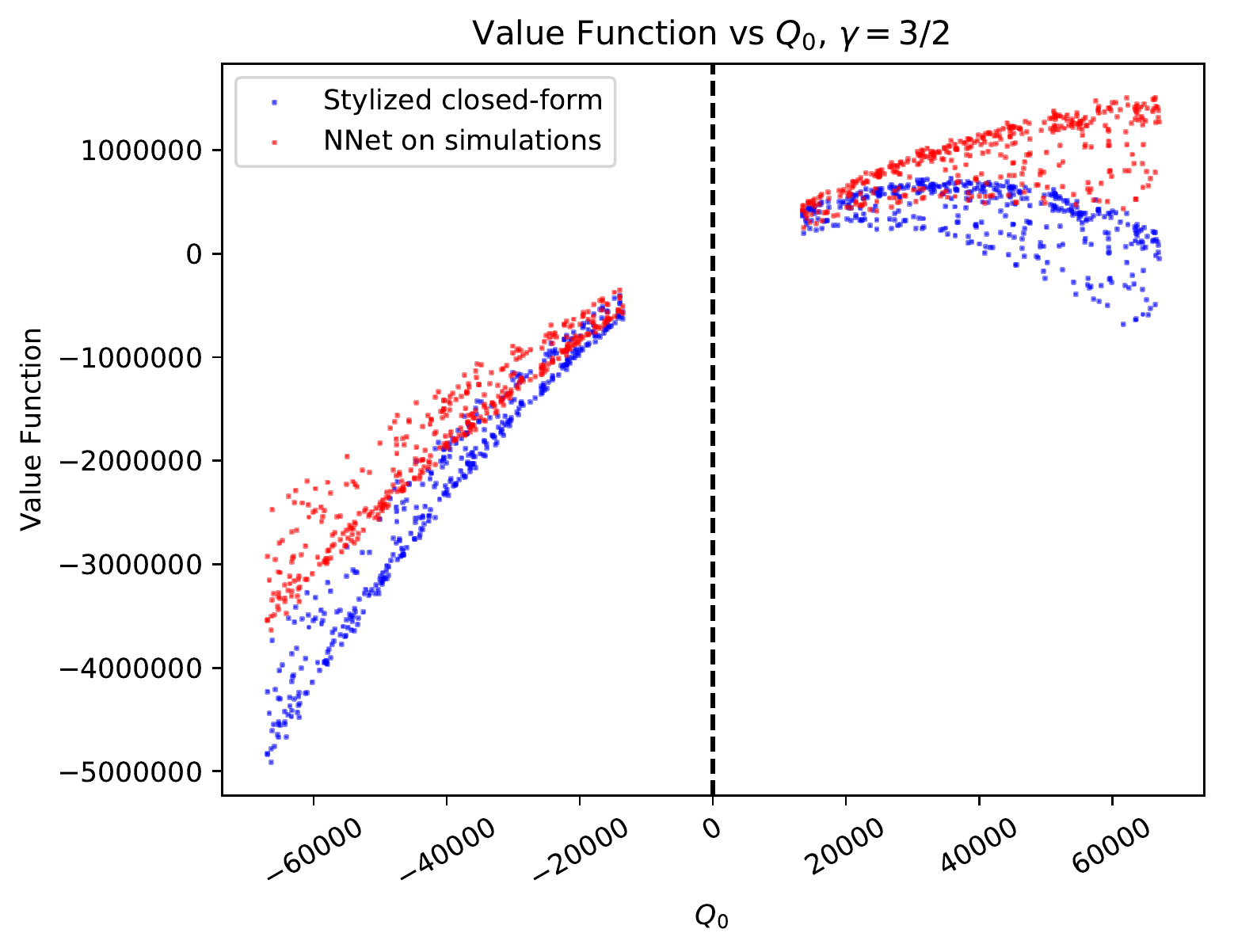}
\end{minipage}
\begin{minipage}{.49\linewidth}
  \includegraphics[width=\linewidth]{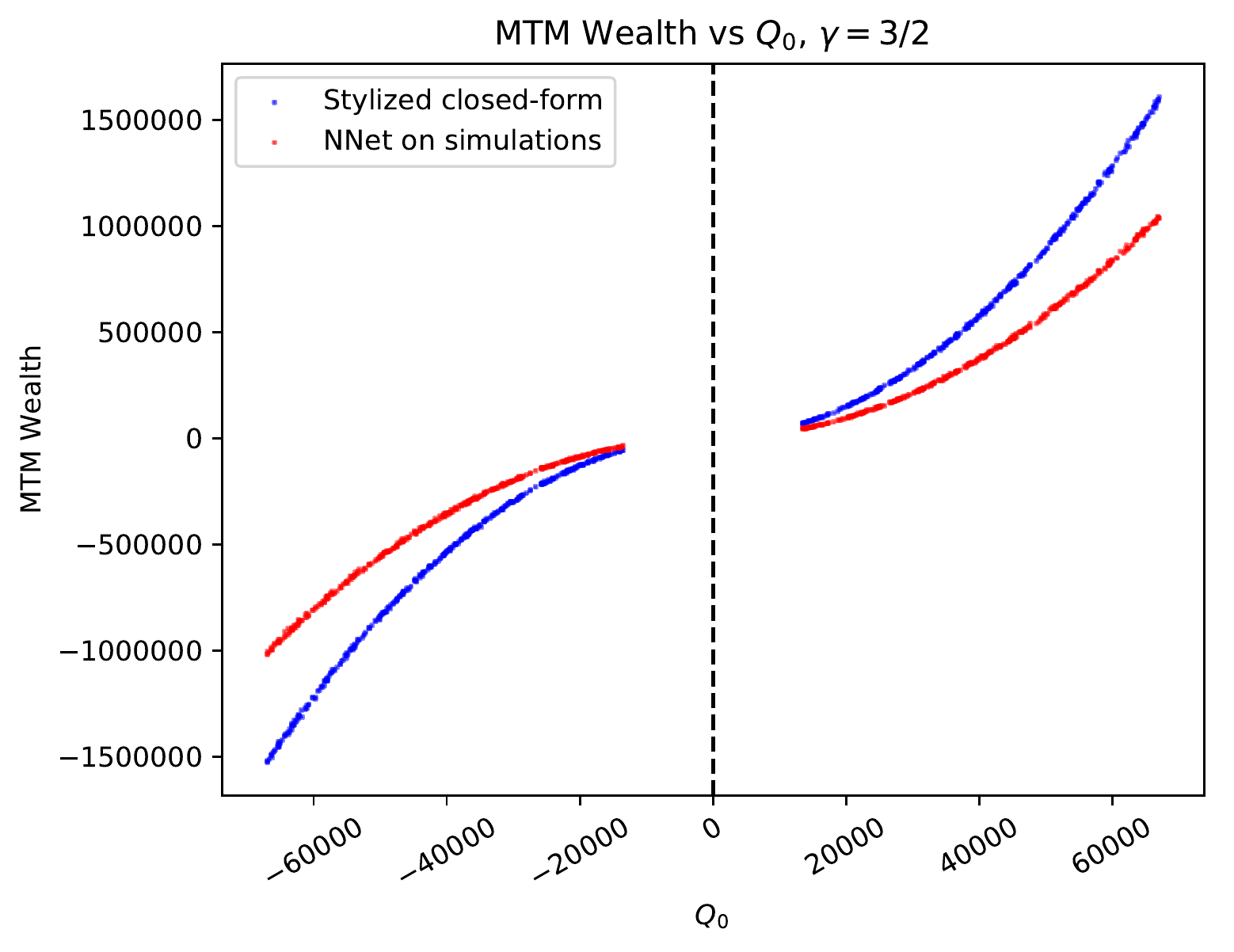}
\end{minipage}
\caption{Value Functions (left) and Marked-to-Market Wealth (right): for $\gamma=2$ (top) and $\gamma=3/2$ (bottom)}
\label{fig:cmp-baseline}
\end{figure}

In Figure~\ref{fig:cmp-baseline}, we compare the performance of the control learnt by the neural network with state-of-the-art control. For the sake of illustration, we focus here on two cases: $\gamma = 2$ and $\gamma = 3/2$.  Given $q_0$ and a choice of control, using Monte Carlo samples we compute an approximation of the reward function (corresponding to~\eqref{eq:1.4} conditioned on $Q_0=q_0$) and the marked-to-market final wealth in cash. In the case $\gamma=2$, using the neural network control (with or without seasonality, and also in the multi-preference setting) almost matches the optimal one obtained by the closed-form solution. When $\gamma=3/2$, the control stemming from the neural network performs better than the control coming from the closed-form solution (recall that the goal is to maximize the reward function). 
We stress that while the reward function is asymmetric due to the intrinsic asymmetry in equation \eqref{eq:1.3}, the marked-to-market wealth defined by $\text{MTM}_T=\text{sign}(Q_0)(Q_0S_0-X_T)$ is symmetric except for a sign adjustment.

\subsection{Empirical Verification of DNN Explainability}

\begin{figure}
\centering
\begin{minipage}{.47\linewidth}
  \includegraphics[width=\linewidth]{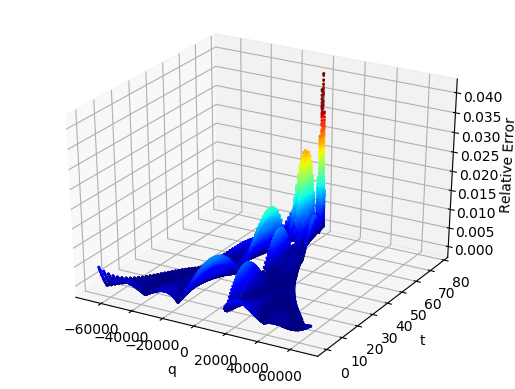}
\end{minipage}
\begin{minipage}{.47\linewidth}
  \includegraphics[width=\linewidth]{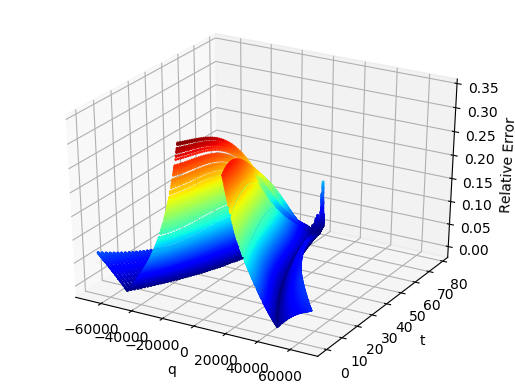}
\end{minipage}
\caption{Relative error in the regression: $\gamma=2$ (left) and $\gamma=3/2$ (right)}
\label{fig:-relerr-DNN}
\end{figure}

Figure~\ref{fig:-relerr-DNN} displays the DNN control and relative error in the neural network regression, namely, if $(t,q) \mapsto f_{\theta}(t,q)$ and $(t,q) \mapsto \tilde f_{\theta}(t,q)$ denote respectively the trained neural network and approximate version obtained by regression, we compute: $(t,q) \mapsto |(f_{\theta}(t,q) - \tilde f_{\theta}(t,q)) / \tilde f_{\theta}(t,q)|$.
The maximum relative error for the case $\gamma=2$ is $0.04$, while for the case $\gamma=3/2$ it is $0.35$. 
We see that for most time steps, the relative error is small in the bulk of the distribution of $q_t$. However, around the terminal time $=78$, the distribution is very concentrated around $q=0$, which explains why the regression is less accurate.

\section{Conclusion}
\label{sec5}
In this paper, we succeed in using a neural network to learn the mapping between end-user preferences and the optimal trading speed, allowing to buy or sell a large number of shares or contracts on financial markets. Prior to this work, various proposals have been made but the learned controls have always been specialized for a given set of preferences. Here, our \emph{multi-preferences neural network} learns the solution of a class of dynamical systems.

Note that optimal execution dynamics are reacting to the applied control, via the price impact of the buying or selling pressure. The neural network hence has to learn this feedback. Our approach uses a deep neural network, whose inputs are user preferences and the state of the optimization problem that change after each decision. The loss function can only be computed at the end of the trading day, once the neural controller has been used 77 times. The backpropagation hence takes place across these 77 steps. We faced some exploration - exploitation issues and solved it by choosing suitable the range of users’ preferences to ensure that long enough trajectories are observed during the learning. 

Our setup leverages on transfer learning, starting on simulated data before switching to historical data. Since we want to understand how the learned controls are different from the closed-form solution on a stylized model that are largely used by practitioners, we learn on different versions of simulated data: from a model corresponding to the stylized one to one incorporating non-stationarities, and then to real data. 

To ensure the explainability of our learned controls, we introduce a projection method on the functional space spanned by the closed-form formula. It allows us to show that most of the learned controls belong to this manifold. When we depart from the traditional model by incorporating a risk aversion term, the learned control is close to its projected version on the bulk of the distribution, supporting the idea that the projection explains a significant part of the neural network control. Furthermore, the source of adaptation to realistic dynamics focuses on a "shift term" $h_1$ in the control. The most noticeable adaptation of $h_1$ exhibits slow oscillations that are most probably reflecting the seasonalities of intra-day price dynamics. 

We then introduce a version of the loss function that reflects more the reality of intra-day price dynamics (that are sub-diffusive). Despite the fact that the associated HJB has no closed form solution, we manage to learn its associated optimal control and show, using our projection technique, that it almost belongs to the same manifold. 
This approach delivers explainability of learned controls and can probably be extended to contexts other than optimal trading. It should help regulators to have a high level of trust in the learned controls, that are often considered as “black boxes”: our proposal exposes the fraction of the controls that belongs to the manifold practitioners and regulators are familiar with, allowing them to perform the usual "stress tests" on it, and quantifies this fraction as a $R^2$ curve that is easy to interpret. More details on the ethical impact and broad societal implications are discussed in Appendix~\ref{sec:ethics}.

\clearpage

\bibliographystyle{rQUF}
\bibliography{my_papers}




\pagebreak 

\appendices
\section{Details on Explicit Solution and Implementation}

\subsection{Details on the explicit solution}

For the sake of completeness, we recall how the benchmark solution is obtained; see~\citep{cj16} for more details. 
The continuous form of the problem defined in equations (1) to (4) of the paper, when $\gamma=2$, can be characterized by the value (we drop the subscripts $A$ and $\phi$ to alleviate the notations)
$$
    \EE_{X_0, S_0, Q_0}[V(0, X_0, S_0, Q_0)],
$$
where $V$ is the value function, defined as:
\begin{equation}
\begin{aligned}
V(t_0,&x,s,q) = \sup_{\nu} \EE \Big[ X_T + Q_T S_T-A |Q_T|^2 -\phi \int_{t_0}^T |Q_t|^2 \Big] \\ 
&\textrm{subject to} \quad \begin{cases}
dS_t = \alpha (\mu_t + \nu_t) dt + \sigma dW_t\\
dQ_t = \nu_t dt\\
dX_t = -\nu_t(S_t+\kappa \nu_t)dt\\
S_t>0, \quad \forall t\\
X_{t_0}=x, Q_{t_0}=q, S_{t_0}=s.
\end{cases}
\end{aligned}
\end{equation}

From dynamic programming, we obtain that the value function $V$ satisfies the following Hamilton-Jacobi-Bellman (HJB): for $t \in [0,T), x,s,q \in \mathbb{R}$,
\begin{align}
    \notag
    \partial _t V &- \phi q^2 + \frac{1}{2} \sigma ^2\partial ^2_S V + \alpha \mu \partial_S V 
    \\
     &+\sup_{\nu}\Big\{\alpha \nu \partial_S V +   \nu \partial_q V -\nu(s+\kappa \nu)\partial _X V \Big\}=0
\label{eq:A.1}
\end{align}
with terminal condition  $V(T,x,s,q) = x+q(s-Aq)$. 

If we use the ansatz $V(t,x,s,q) = x+qs+u(t,q)$, with $u$ of the form $u(t,q) = h_0(t) + h_1(t)q + h_2(t)\frac{q^2}{2}$, the optimal control resulting from solving this problem can be written as:
\begin{equation}
\nu^* (t,q) = \frac{\alpha q + \partial _q u(t,q)}{2 \kappa}
= \frac{h_1(t)}{2 \kappa}  + \frac{\alpha +h_2(t)}{2 \kappa} q.
\end{equation}
Hence, $h_1(t)$ and $h_2(t)$ act over the control by influencing either the intercept of an affine function of the inventory, in the case of $h_1(t)$, or its slope, in the case of $h_2(t)$. 

From the HJB equation, these coefficients are characterized by the following system of ordinary differential equations (ODEs):
\begin{align}
\begin{cases}
\dot{h_2}(t) = \big(2 \phi - \frac{1}{2\kappa}\alpha ^2 \big) - \frac{\alpha}{\kappa} h_2(t) - \frac{1}{2\kappa}h_2^2(t),  \\
\dot{h_1}(t) +\frac{1}{2\kappa}(\alpha+ h_2(t))h_1(t) = -\alpha \mu(t), \\
\dot{h_0}(t) = - \frac{1}{4\kappa}h_1^2(t)
\label{eq:A.15}
\end{cases}
\end{align}

with terminal conditions:
\begin{align}
\begin{cases}
h_0(T) = 0,\\
h_1(T) = 0,\\
h_2(T) = -2A.
\end{cases}
\label{eq:A.19}
\end{align}

\clearpage
\subsection{Details on the implementation}
 
In this section, we provide more details on the implementation of the method based on neural network approximation. For the neural network, we used a fully connected architecture, with three hidden layers, five nodes each.

\begin{figure}
\centering
\resizebox{.8\textwidth}{!}{%
\begin{tikzpicture}[
roundnode/.style={circle, draw=red!60, fill=blue!5, very thick, minimum size=10mm},
squarenode/.style={rectangle, draw=blue!60, fill=blue!5, very thick, minimum size=6mm},
typetag/.style={rectangle, draw=cyan!50, ultra thin, font=\scriptsize\ttfamily, anchor=west},
typetag2/.style={rectangle, draw=blue!50, very thick, font=\scriptsize\ttfamily, anchor=west}
]

\node[roundnode](ctrl){$\nu_0$};

\node[squarenode](t0)[above left=11mm and 1cm of ctrl.west]{$t_0$};
\node[squarenode](Q0)[above left=1mm and 1cm of ctrl.west]{$Q_{t_0}$};
\node[squarenode](X0)[below left=1mm and 1cm of ctrl.west]{$X_{t_0}$};
\node[squarenode](S0)[below left=11mm and 1cm of ctrl.west]{$S_{t_0}$};

\node (fit_a0)[draw=black!50,typetag, fit={(t0) (Q0) (X0) (S0)}] {};
\node (fit_b0)[draw=black!50,typetag2, fit={(t0) (Q0)}] {};

\node[squarenode](t1)[above right=11mm and 1cm of ctrl.east]{$t_1$};
\node[squarenode](Q1)[above right=1mm and 1cm of ctrl.east]{$Q_{t_0}$};
\node[squarenode](X1)[below right=1mm and 1cm of ctrl.east]{$X_{t_0}$};
\node[squarenode](S1)[below right=11mm and 1cm of ctrl.east]{$S_{t_0}$};
  
\node (fit_a1)[draw=black!50,typetag, fit={(t1) (Q1) (X1) (S1)}] {};
\node (fit_b1)[draw=black!50,typetag2, fit={(t1) (Q1)}] {};

\node[roundnode](ctrl_1)[right=4cm of ctrl.west]{$\nu_1$};

\node[squarenode](t2)[above right=11mm and 1cm of ctrl_1.east]{$t_2$};
\node[squarenode](Q2)[above right=1mm and 1cm of ctrl_1.east]{$Q_{t_2}$};
\node[squarenode](X2)[below right=1mm and 1cm of ctrl_1.east]{$X_{t_2}$};
\node[squarenode](S2)[below right=11mm and 1cm of ctrl_1.east]{$S_{t_2}$};
  
\node (fit_a2)[draw=black!50,typetag, fit={(t2) (Q2) (X2) (S2)}] {};
\node (fit_b2)[draw=black!50,typetag2, fit={(t2) (Q2)}] {};

\node(dots) [right=of fit_a2] {$\cdots$};

\node[roundnode](ctrl_last)[right=10cm of ctrl.west]{$\nu_{T-1}$};

\node[squarenode](T)[above right=11mm and 1cm of ctrl_last.east]{$T$};
\node[squarenode](QT)[above right=1mm and 1cm of ctrl_last.east]{$Q_T$};
\node[squarenode](XT)[below right=1mm and 1cm of ctrl_last.east]{$X_T$};
\node[squarenode](ST)[below right=11mm and 1cm of ctrl_last.east]{$S_T$};
  
\node (fit_aT)[draw=black!50,typetag, fit={(T) (QT) (XT) (ST)}] {};
\node (fit_bT)[draw=black!50,typetag2, fit={(T) (QT)}] {};

\draw[->,typetag2](fit_b0.east) -- (ctrl.west);
\draw[->,typetag](ctrl.east) |- (fit_a1.west);

\draw[->,typetag2](fit_b1.east) -- (ctrl_1.west);
\draw[->,typetag](ctrl_1.east) |- (fit_a2.west);

\draw[->,typetag2](fit_b2.east) -- (dots.west);
\draw[->,typetag](ctrl_last.east) |- (fit_aT.west);

\end{tikzpicture}
}
\caption{Structure of the State Variables' Simulation - All Steps}
\label{fig:4}
\end{figure}
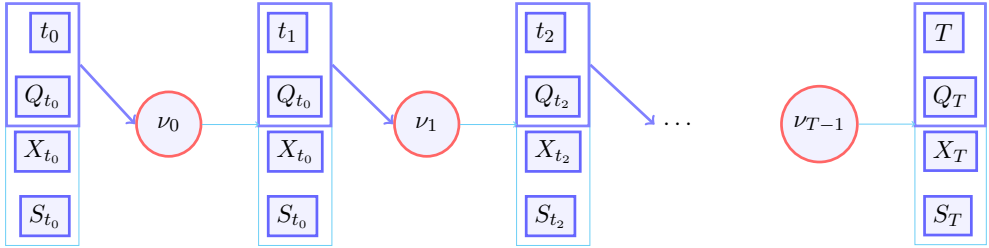

One forward step in this setup is described by Figure \ref{fig:3}, while one round of the SGD is represented by Figure \ref{fig:4}. The \emph{same} neural network learns from the state variables obtained in the previous step. It thus learns the influence of its own output through many time steps. We have experimented with using one neural network per time step. However, this method implies more memory usage, and did not provide any clear advantage with respect to the one presented here.

From Figures \ref{fig:3} and \ref{fig:4}, we can clearly see the idea that the control influences the dynamics. We are, in fact, optimizing in a closed-loop learning environment, where the trader's actions have both permanent and temporary market impact on the price.

The mini-batch size we used, namely 64, is relatively small. While papers like \citep{dcmc12} defend the use of larger mini-batch sizes to take advantage of parallelism, smaller mini-batch sizes such as the ones we use are shown to improve accuracy in \citep{kmn16}, \citep{ml18} and \citep{wm03}. 

\begin{figure}
\centering
\resizebox{.3\textwidth}{!}{ 
\begin{tikzpicture}[
roundnode/.style={circle, draw=green!60, fill=blue!5, very thick, minimum size=7mm},
roundnode2/.style={circle, draw=red!60, fill=blue!5, very thick, minimum size=10mm},
squarenode/.style={rectangle, draw=blue!60, fill=blue!5, very thick, minimum size=8mm},
]

\node[roundnode2](ctrl){$\nu_t$};
\node[squarenode](input_t)[above left=6mm and 2cm of ctrl.east]{$t_i$};
\node[squarenode](input_Q)[below left=6mm and 2cm of ctrl.east]{$Q_{t_i}$};

\node[squarenode](output_Q)[above right=15mm and 3cm of ctrl.west]{$Q_{t_{i+1}}$}; 
\node[squarenode](output_X)[right=3cm of ctrl.west]{$X_{t_{i+1}}$};
\node[squarenode](output_S)[below right=15mm and 3cm of ctrl.west]{$S_{t_{i+1}}$};

\node[roundnode](kappa) at ($(output_Q)!0.5!(output_X)$)[left=1cm] {$\kappa_t$};
\node[roundnode](alpha) at ($(output_X)!0.5!(output_S)$)[left=1cm] {$\alpha_t$};
\node[roundnode](epsilon)[below right=26mm and 18mm of ctrl.west]{$\epsilon_t$}; 

\draw[-stealth', rounded corners=5mm](input_t.south) |- (ctrl.west);
\draw[-stealth', rounded corners=5mm](input_Q.north) |- (ctrl.west);

\draw[-stealth', rounded corners=5mm](kappa.south) |- (output_X.west);
\draw[-stealth', rounded corners=5mm](ctrl.east) -- (output_X.west);

\draw[-stealth', rounded corners=5mm](alpha.south) |- (output_S.west);
\draw[-stealth', rounded corners=5mm](epsilon.north) |- (output_S.west);
\draw[-stealth', rounded corners=5mm](ctrl.south) |- (output_S.west);

\draw[-stealth', rounded corners=5mm](ctrl.north) |- (output_Q.west);
\end{tikzpicture}
}
\caption{Structure of the state variables' simulation - One Step}
\label{fig:3}
\end{figure}
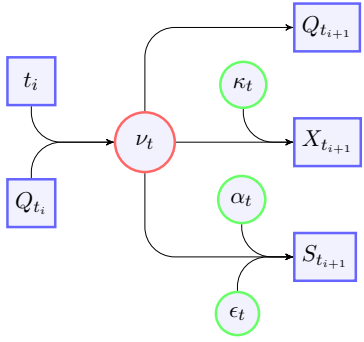

Simulations are run on a MacOS Mojave laptop with 2.5 GHz Intel Core i7 and 16G of RAM, without GPU acceleration. Available GPU cluster did not increase the average speed of the simulations. Tests done on CPU clusters composed of Dell 2.4 GHz Skylake nodes also did not indicate relevant speed improvements.

\subsection{Kurtosis and Auto-correlation of stock returns}

Table \ref{tab:2.3} shows the values of kurtosis and auto-correlation for all the stocks. They respectively indicate that the returns have heavy tails and are auto-correlated. Neither of these characteristics of real stock returns are taken into account in the baseline model, but they are easily accounted for in the deep neural network setup presented in the paper.

\begin{table}
    \centering
    \tbl{Kurtosis and auto-correlation of intra-day returns}{%
        \begin{tabular}{|l|l|l|l|l|l|l|l|}
        \hline
        Kurtosis            & ABX   & AEM   & AGU  & BB    & BMO   & BNS  & COS  \\
                            & 222   & 264   & 1825 & 2344  & 1796  & 2087 & 2528  \\\cline{2-8}
                            & DOL   & GIL   & GWO   & HSE  & MRU   & PPL  & RCI.B \\
                            & 163   & 7746  & 507   & 2975 & 319   & 7006 & 2188 \\\cline{2-8}
                            & RY    & SLF   & TRI  & TRP   & VRX  && \\ 
                            & 2648  & 11302 & 403  & 7074  & 213  && \\ 
        \hline\hline
        Auto-correlation    & ABX   & AEM   & AGU  & BB    & BMO   & BNS  & COS  \\
        (lag = 5 min)       & -.025 & -.026 & -.03 & -.031 & -.026 & -.03 & -.041\\\cline{2-8}
                            & DOL   & GIL   & GWO   & HSE  & MRU   & PPL  & RCI.B \\
                            & -.118 & -.047 & -.094 & -.02 & -.105 & -.15 & -.062\\\cline{2-8}
                            & RY    & SLF   & TRI  & TRP   & VRX  && \\  
                            & -.017 & -.07  & -.04 & -.066 & -.075 &&\\ \hline
        \end{tabular}
     }%
\label{tab:2.3}
\end{table}

\clearpage
\section{Learning the mapping to the optimal control in a closed loop}
\label{sec4}

\subsection{Details on market impact parameters estimation}

\paragraph{Permanent market impact}
Let $\mathcal{S}$ be the space of stocks; and $\mathcal{D}$ be the space of trading days. If we take five-minute intervals (indicated by the superscript notation), we can write equation (1) for each stock $s \in \mathcal{S}$, for each day $d \in \mathcal{D}$, and for each five-minute bin indexed by $t$ as:
\begin{align}
\Delta S_{s,d,t}^{\text{5min}} = \alpha_{s,t} \mu_{s,d,t}^{\text{5min}} \Delta t + \sigma_s \sqrt{\Delta t} \epsilon_{s,d,t}^{\text{5min}},
\label{eq:2.3}
\end{align}
\noindent where the subscripts $s,d,t$ respectively indicate the stock, the date, and the  five-minute interval to which the variables refer, and $\Delta t = $ 5 min, by construction. We have $\alpha_{s,t}$ independent of $d$, which assumes that for any given day the permanent market impact multiplier the agent may have on the price of a particular stock $s$, for a given time bin $t$ is the same. And we have $\sigma_s$ independent of the day $d$ and time bin $t$, which means the volatility related to the noise term is constant for a given stock. 

As an empirical proxy for the theoretical value $\mu_{s,d,t}^{\text{5min}}$ representing the net order flow of all the agents, we use the order flow imbalance observed in the data: $\text{Imb}_{s,d,t}^{\text{5min}}$. We define this quantity as:
\begin{align}
    \text{Imb}_{s,d,t}^{\text{5min}} = \sum_{j=t}^{t+\text{5min}} v^{\text{buy}}_{s,d,j}  - \sum_{j=t}^{t+\text{5min}} v^{\text{sell}}_{s,d,j},
    \label{eq:2.4}
\end{align}

\noindent where $v^{\text{buy}}_{s,d,t}$ is the volume of a buy trade, and $v^{\text{sell}}_{s,d,t}$ is the volume of a sell trade, and we aggregate their net amounts over five minute bins.

However, since we estimate the permanent market impact parameter $\alpha$ using data from different stocks, and we would like to perform a single regression for all of them, we re-scale the data used in the estimation. In order to make the data from different stocks comparable, we do the following:

\begin{enumerate}
\itemsep-0.3em
    \item Divide the trade price difference ($\Delta S_{s,d,t}^{\text{5min}}$) by the average bin spread over all days for the given stock, also calculated in its respective five minute bin:
    \begin{align}
        \overline{\Delta S}^{\text{5min}}_{s,d,t} = \frac{\Delta S_{s,d,t}^{\text{5min}} }{\frac{1}{|\mathcal{D}|}\sum_{d\in \mathcal{D}} \psi_{s,d,t}^{\text{5min}}},
        \label{eq:2.5}
    \end{align}
    where $\psi_{s,d,t}^{\text{5min}}$ is the bid-ask spread for a given (stock, date, bin) tuple. 
    
    \item Divide the trade imbalance ($\text{Imb}_{s,d,t}^{\text{5min}}$) by the average bin volume over all days for the given stock, calculated on the respective five minute bin:
    \begin{align}
        \overline{\text{Imb}}^{\text{5min}}_{s,d,t} = \frac{\text{Imb}_{s,d,t}^{\text{5min}}}{\frac{1}{|\mathcal{D}|}\sum_{d\in\mathcal{D}}(\text{Total Volume}_{s,d,t}^{\text{5min}})}.
        \label{eq:2.7}
    \end{align}
    where $\text{Total Volume}_{s,d,t}^{\text{5min}}$ stands for the total traded volume for both buys and sells for a given (stock, date, bin) tuple.
\end{enumerate}

This ensures that both the volume and the price are normalized to quantities that can be compared. Each $(s,d,t)$ now gives rise to one data point, and we can thus run the following regression instead, using comparable variables, to find $\bar{\alpha}$: 
\begin{align}
\overline{\Delta S}^{\text{5min}}_{s,d,t} = \bar{\alpha} \cdot \overline{\text{Imb}}^{\text{5min}}_{s,d,t} + \bar{\epsilon}_{s,d,t}^{\text{5min}},
\label{eq:2.9}
\end{align}

\noindent where $\bar{\alpha}$ is the new slope parameter we would like to estimate, and $\bar{\epsilon}_{s,d,t}^{\text{5min}}$ is the normalized version of the residual we had in equation \eqref{eq:2.3}. 

In order to use $\bar{\alpha}$ in a realistic way, we de-normalize the regression equation for each stock by doing:
\begin{align}
    \Delta S^{\text{5min}}_{s,d,t} = \Bigg(\frac{\bar{\alpha}}{\Delta t} \frac{\bar{\psi}_{s,t}}{\bar{\mathcal{V}}_{s,t}}\Bigg) \cdot  \text{Imb}^{\text{5min}}_{s,d,t}  + \sigma \sqrt{\Delta t} \epsilon_{s,d,t}^{\text{5min}}.
    \label{eq:2.11}
\end{align}

\noindent where $\bar{\mathcal{V}}_{s,t}$ is the average bin volume for stock $s$ at bin the time bin $t$, $\bar{\psi}_{s,t}$ is the average bin bid-ask spread for the same pair $(s,t)$.

The derivation for the \textbf{temporary market impact} $\kappa$ uses equation (2) in the paper, and follows similar steps.

\paragraph{Temporary market impact} The $\kappa$ parameter represents the magnitude of the trading cost of an agent in the market. This is an indirect cost incurred by the trader due to the impact that their order has in the market mid-price. From the dynamics of the wealth process $dX_t = -\nu_t(S_t+\kappa \nu_t)dt, \label{eq:2.12}$
we are assuming that the trading cost is linear in the speed of trading. Notice that the faster the agent needs to trade, the larger the cost incurred, and the opposite is true for slower trading. We can re-normalize \eqref{eq:2.12} by rewriting it as:
\begin{align}
     \delta x := (-1)\cdot \frac{dX_t}{\nu dt} = S_t + \kappa \nu_t.
    \label{eq:2.13}
\end{align}

Then, for each stock $s  \in \mathcal{S}$, for each day $d \in \mathcal{D}$, and for each five-minute bin indexed by $t$, we would like to estimate the following expression:
\begin{align}
    \delta x_{s,d,t}^{\text{5min}} = S_{s,d,t,\text{start}}^{\text{5min}} + \kappa \nu_{s,d,t}^{\text{5min}} + \epsilon_{s,d,t}^{\text{5min}}.
    \label{eq:2.14}
\end{align}
\noindent where the subscript `start' indicates that we are using the mid-price at the start of the five-minute bin indexed by $t$. 

Since we do not have access to $\nu$, as it is private information for each trader, we would like to compensate by restricting our statistical analysis to bins that are dominated by only one trader, but not so much as to have stagnant prices. With the goal of isolating one trader from the crowd in mind, we seek to restrict our dataset based on the imbalance and on the size of the dominating traded volume in a given bin. 

We extract the trade volume for each bin as the maximum between the aggregate buy volume and the aggregate sell volume for that particular bin. With this criteria we simultaneously define whether we have a buy or a sell, and thus the sign of the trade (+1 for buys, and -1 for sells). In mathematical terms, the dominant trading volume by bin, $D_{s,d,t}^{\text{5min}}$, can be expressed as:
\begin{align}
    D_{s,d,t}^{\text{5min}} = \frac{\max_{\text{5min}}\Big\{\sum_{j=t}^{t+\text{5min}} v^{\text{buy}}_{s,d,j},  \sum_{j=t}^{t+\text{5min}} v^{\text{sell}}_{s,d,j}\Big\}}{{\frac{1}{|d|}\sum_{d \in \mathcal{D}}(v_{s,d,t}^{\text{5min}})}}.
    \label{eq:2.16}
\end{align}
\noindent where all variables have been defined in the paragraph about Permanent market impact. 

As \cite{pb18} find, aggregate-volume impact saturates for large imbalances, on all time scales, and highly biased order flows are associated with very small price changes. They find that the probability for an order to change the price decreases with the local imbalance, and vanishes when the order signs are locally strongly biased in one direction. However, if we first restrict the data to having the dominant trading volume by bin $D_{s,d,t}^{\text{5min}}$, the imbalances are naturally restricted. Hence, we do not need to further restrict the upper bound of the imbalance. However, we do want to restrict its lower bound, since we want to analyse bins that isolate either the buying or the selling activity of only one trader. 

In short, to estimate $\kappa$ we restrict our dataset in the following way:
\begin{enumerate}
    \item $2<D_{s,d,t}^{\text{5min}} <10$,
    \item $|\text{Imb}_{s,d,t}^{\text{5min}}| > |\text{Imb}_{s,d,t}^{\text{5min}}| _{40\%}$.
\end{enumerate}

For the trading cost $\kappa$ estimation, we would like to track agents who traded more than usual during five minutes. Hence, if the volume traded in the direction that had most trades is not minimally significant, we would not like to use it in the estimation. We have defined in equation \eqref{eq:2.16} the dominant trading volume as being the volume in the direction that had more trades during a given bin (either buy or sell). 

We exclude from the sample those trades that had a dominant trading volume lower than twice the average dominant volume for a given bin. Since we only want to keep track of one individual agent, we also restrict this same quantity to be below ten times the average dominant trade volume for the bin. We believe that if the volume is larger than this, then it is very likely that other agents are also involved in the trade. 

The third and last filter we make is for the trade imbalance, as we would like to isolate the effect of buys and sells. We keep in the sample only those observations on which the absolute value of the imbalance is greater than the 40\% percentile of the signed imbalances.

Once these filters have been put into place, we are almost ready to estimate the agent's transaction costs. For our new subset of  $(s,d,t)$ only, we would like to establish:
\begin{align}
    \delta x_{s,d,t}^{\text{5min}} - S_{s,d,t,\text{start}}^{\text{5min}} = \kappa \nu_{s,d,t}^{\text{5min}}.
    \label{eq:2.17}
\end{align}

Nevertheless, we would like to do it for all stocks at the same time, just like we did for the parameter $\alpha$. Hence, we need to re-normalize both sides of this equation to comparable values. For the left-hand side, we re-normalize price moves by the average bid-ask spread of a given stock on a given bin. For the right-hand side, we re-normalize $\nu_{s,d,t}^{\text{5min}}$ by the average trading volume at the bin. Notice that this is precisely $D_{s,d,t}^{\text{5min}}$, as defined in equation \eqref{eq:2.16}. Hence, we have:
\begin{align}
    \frac{\delta x_{s,d,t}^{\text{5min}} - S_{s,d,t,\text{start}}^{\text{5min}}}{\frac{1}{d}\sum_{d \in \mathcal{D}} \psi_{s,d,t}^{\text{5min}}} = \kappa \frac{\nu_{s,d,t}^{\text{5min}}}{\frac{1}{|d|}\sum_{d \in \mathcal{D}} v_{s,d,t}^{\text{5min}}}.
    \label{eq:2.18}
\end{align}

Having re-normalized our variables, we can now run the regression on equation \eqref{eq:2.12} to find $\kappa$ and estimate individual trading costs. In practice, we will be using the VWAP defined in section (2.1), item 4 as the input for $\delta x$. 

The regression we estimate in order to find $\kappa$ is, thus:
\begin{align}
    \frac{\text{sign}_{s,d,t}^{\text{5min}} \cdot \large(\text{VWAP}_{s,d,t}^{\text{5min}}-S_{s,d,t,\text{start}}^{\text{5min}}\large)}{\frac{1}{|d|}\sum_{d \in \mathcal{D}} \psi_{s,d,t}^{\text{5min}}} = \kappa \frac{\nu_{s,d,t}^{\text{5min}}}{\frac{1}{|d|}\sum_{d \in \mathcal{D}} v_{s,d,t}^{\text{5min}}}.
    \label{eq:2.19}
\end{align}

\subsection{Details on the multi-preferences controller}

In section 2.3 of the paper, we discussed the exploration - exploitation trade-off we encountered when projecting the neural network controller onto the manifold of closed-form solutions when $\gamma=2$.

\begin{figure}
\centering
\begin{minipage}{.49\linewidth}
  \includegraphics[width=\linewidth]{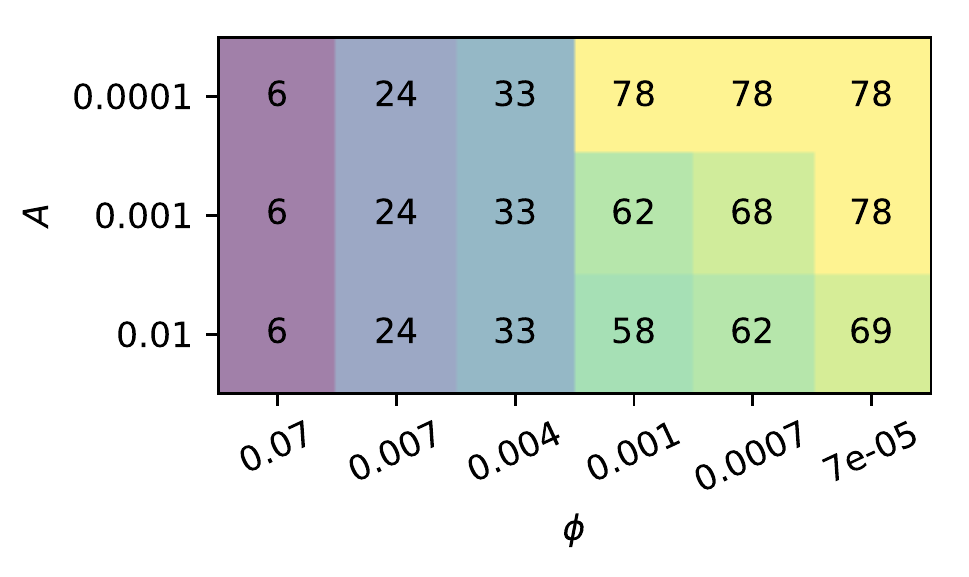}
\end{minipage}
\begin{minipage}{.49\linewidth}
  \includegraphics[width=\linewidth]{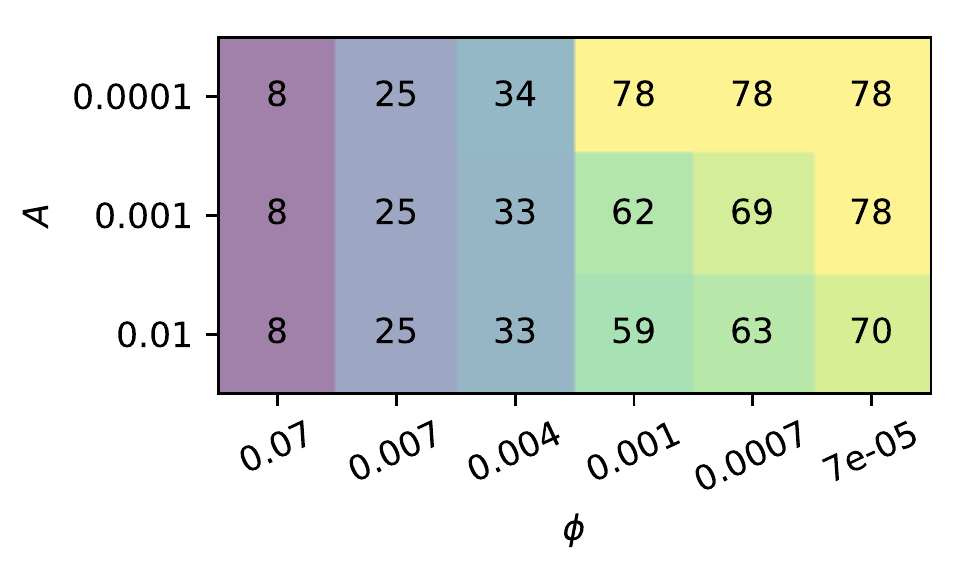}
\end{minipage}
\caption{Time steps until at least 90\% of the order is executed}
\label{fig:5}
\end{figure}

In Figure \ref{fig:5}, we compare the number of time steps it takes the controller to trade at least 90\% of the initial inventory in the market. Rows stand for values of the risk parameter $A\in\{0.01,0.001,0.0001\}$, while columns represent the risk parameter $\phi \in \{0.07,0.007,0.004,0.001,7e^{-4},7e^{-5}\}$.

When the trader's preferences over risk aversion gives them incentive to trade faster, the inventory is executed sooner. The controller is `exploiting' the available trading choices, and the projection on the closed-form manifold becomes harder to estimate. On the other hand, if the pair $(A,\phi)$ defines less restrictive preferences towards execution, there will be some (although not a lot of) inventory left at the end of the trading day that allows us to learn the projection accurately.

\clearpage

\section{Ethical impact and broad societal implications}
\label{sec:ethics}

\paragraph{Explainability of learned control for financial markets.}
Attempts of using machine learning for financial markets are booming the last 5 years. It spans a wide spectrum of applications: from client profiling to nowcasting using databases of texts and satellite images. One area where the acceptance of ML goes slowly is the automation of hedging strategies. These strategies are part of the functioning of financial markets in the post 2008-crisis environment due to genuine regulations demand market participants to compute and compensate their exposure to as many risk scenarios as possible. Moreover, and for obvious reasons, regulators ask for a better understanding of algorithms that look like "black boxes" before allowing the use of ML in production to improve the efficiency of hedging. It is related to the very well-known topics of ``\emph{explainability of AI}" and ``\emph{ethics of algorithms}". 

Optimal trading addresses part of these hedging strategies. A framework like the one in this paper corresponds to a ``large investor", typically a pension fund, who decided to rebalance its portfolio for good reasons (think about the market moves following the spread of COVID-19: the optimal response of pension funds is clearly to rebalance their portfolio to get less exposure to factors like airplane companies, to prevent their pensioners to suffer from losses in the coming year), but they have to do it at a slow enough pace to not push the price too much at their own detriment \citep{ll18}. The reason for the big COVID-19 drop of markets followed by a rebound 10 days later is due to non-optimal trading strategies: the institutions sold too fast, paying too much for their deleveraging and perturbing the prices with the downward pressure of accelerated sells. Being able to adapt the trading strategies to realistic features of intra-day price dynamics (like seasonality and auto-correlations) decreases both the trading costs for large institutions like pension funds (and consequently to their pensioners) and is profitable to the public price formation, that is otherwise does not reflect the ``fair value" of traded instruments.

In the introduction, we cite papers on improving optimal trading, via reinforcement learning or directly using deep learning. They do not attempt to \emph{provide an explainable version of the obtained controls}. In this paper we do it in a very practical way that is compliant with the practices of the financial industry: we project the learned controls on the manifold spanning the controls currently in use by brokers and asset managers, and measure the distance between the obtained projection and the learned controls. In doing this, we provide evidences that the learned controls that we generated are largely contained in this regular manifold, and we allow practitioners and regulators to apply their usual stress testing and certification practices to the projected controls. It is the first proposal in this direction, and we hope that it will propagate in other areas of finance, like deep hedging. It may also be used in other industries when the regulatory demand in explainability is high.

\paragraph{``Functional learning" for control}
The second ``broad impact" of this paper is what we called ``functional learning" in the scope optimal trading with the setup we call \emph{multi-preferences neural network}. In other applications of optimal control, the loss function is fixed because it correspond to one universal use of the controlled object. On financial markets the loss function is parametrized by hyperparameters describing the risk aversion of the end-user. In game theory versions of the optimal trading framework, they are usually called ``agent preferences" \citep{cl18}. 
It is like allowing the driver of an autonomous car to choose the ``driving style" they want to use: far below the speed limit (may be for low carbon emission reasons) or just below, or a ``sport mode", etc. In financial markets, mainly because on the one hand models of price impact have confidence intervals and on the other hand asset owners can have a lot of reasons to rebalance their portfolio (think about an Australian pension fund that has to sell fast enough to face the exceptional pensioners' redemption authorized by the government following the COVID-19 crisis,  see~\cite{Burgess2020}; it is clear that the selling speed has to be fast enough to give the money to pensioners), these parameters are chosen at the start of trading ($t=0$ in our framework) by the end-user. Today a stylized version of the control problem is numerically solved on the fly (either in closed-form, either by a numerical scheme) for the chosen preferences, but the dynamics have to be stylized, hence simplified, enough to allow to solve this very quickly. Learning the optimal control each time an end-user choose new preferences is not fast enough. 

In this paper, the neural net is learning to solve the optimal control problem for any preference; in fact it learns the mapping between the preferences and the optimal solution. It means more flexibility to adapt to user's needs that ``one learning fits all model". 

It effectively learns to solve a whole parametrized family of Hamilton-Jacobi-Bellman equations. Hence it can provide the optimal trading strategy to a given choice of preferences in a flash. Keep in mind this optimal strategy is not one number but it is itself a mapping between the state space of the control problem and the optimal trading speed to be applied.

It had chances to work because it is straightforward to check (in the continuous setting) that the infinite-dimensional optimal control strategy is a smooth function of the two-dimensional vector of preferences $(A,\phi)$, hence we ask to the neural network to interpolate high dimensional functions in this 2-dimensional grid. Nevertheless it is far from easy (that's probably why other attempts did not succeed by now).

In the paper, we refer to what we suspect to be the main difficulty under the name of     ``exploration - exploitation issue". 
For some pairs of user preferences, the optimal trading speed is so fast that the full order is bought or sold far before $T$, and hence the neural network cannot learn from the feedback of its control on the dynamics, simply because it stopped to interact with them.
To counter this we used the closed-form formula of the stylized problem to define a domain of user preferences that is compatible with our $T$, it helped the convergence of the ``functional learning".

\end{document}